\documentclass[10pt,a4paper,twoside]{article}
\usepackage{amsmath, exscale, epsfig,amssymb,makeidx, mathrsfs}
\small
\normalsize
\usepackage[latin1]{inputenc}
\usepackage[T1]{fontenc}

\usepackage{amsfonts}

\usepackage{amsbsy}

\usepackage{amscd}
\usepackage{theorem}

\usepackage{epsf}

\usepackage{psfrag}
\usepackage{graphicx}

\usepackage{multicol}
\allowdisplaybreaks

\def\build#1_#2^#3{\mathrel{\mathop{\kern 0pt#1}\limits_{#2}^{#3}}}
\def\noi{{\noindent}}

\def\un{{\bf 1}}
\newcommand{\bbD}{\mathbb{D}}

\newcommand{\bE}{{\bf E}}

\newcommand{\bm}{{\bf m}}

\newcommand{\bN}{\mathbb{N}}

\newcommand{\bP}{{\bf P}}
\newcommand{\bR}{\mathbb{R}}

\newcommand{\cI}{\mathcal{I}}

\newcommand{\cF}{\mathcal{F}}
\newcommand{\cG}{\mathcal{G}}

\newcommand{\cL}{\mathcal{L}}
\newcommand{\cM}{\mathcal{M}}
\newcommand{\cN}{\mathcal{N}}

\newcommand{\cT}{\mathcal{T}}

\def\baN{{\bf N}}

\def\cJ{{\cal J}}
\def\cD{{\cal D}}

\def\varep{\varepsilon}

\def\be{\begin{equation}}
\def\ee{\end{equation}}
\def\ba{\begin{eqnarray*}}
\def\ea{\end{eqnarray*}}

\def\noi{\noindent}

\newcommand{\lgeo}{[{\!  } [}
\newcommand{\rgeo}{] {\!  } ]}

\def\cqfd{ \hfill $\blacksquare$ }

\newtheorem{theorem}{Theorem}[section]
\newtheorem{lemma}[theorem]{Lemma}
\newtheorem{proposition}[theorem]{Proposition}
\newtheorem{corollary}[theorem]{Corollary}

{\theorembodyfont{\rmfamily}}
{\theorembodyfont{\rmfamily}}
{\theorembodyfont{\rmfamily}}
{\theorembodyfont{\rmfamily}}

\begin{document}

\title{ {\bf  EXCEPTIONALLY SMALL BALLS IN STABLE TREES}\thanks{ This research was supported by the grant ANR A3, Projet BLAN-****}}
\author{Thomas {\sc Duquesne}\thanks{LPMA;
Universit\'e P. et M. Curie (Paris 6), Bo\^ite courrier 188, 4 place Jussieu, 75252 Paris Cedex 05, FRANCE. Email: thomas.duquesne@upmc.fr} \and Guanying {\sc Wang}
\thanks{Department of Mathematical Sciences, Tsinghua University, Beijing, 100084, PR CHINA. Email: wanggy05@tsinghua.edu.cn}}
\vspace{2mm}

\date{\today}

\maketitle

\begin{abstract} The $\gamma$-stable trees are random measured compact metric spaces that appear as the scaling limit of Galton-Watson trees whose offspring distribution lies in a $\gamma$-stable domain, $\gamma \in (1, 2]$. They form a specific class of L\'evy trees (introduced by Le Gall and Le Jan in \cite{LGLJ1}) and the Brownian case $\gamma= 2$ corresponds to Aldous Continuum Random Tree (CRT).
 In this paper, we study fine properties of the mass measure, that is the natural measure on $\gamma$-stable trees. We first discuss the minimum of the mass measure of balls with radius $r$ and we show that this quantity is of order $r^{\frac{\gamma}{\gamma-1}} (\log1/r)^{-\frac{1}{\gamma-1}}$.
We think that no similar result holds true for the maximum of the mass measure of balls with radius $r$,
except in the Brownian case: when $\gamma = 2$, we prove that this quantity is of order $r^2 \log 1/r$.
In addition, we compute the exact constant for the lower local density
of the mass measure (and the upper one for the CRT), which continues previous results from \cite{Du9, Du10, DuLG3}.

\medskip

\end{abstract}

\section{Introduction}
\label{introsec}
Stable trees are particular instances of L\'evy trees that form a class of random compact metric spaces introduced by Le Gall and Le Jan in \cite{LGLJ1} as the genealogy of Continuous State Branching Processes (CSBP for short). The class of stable trees contains Aldous's continuum random tree that corresponds to the Brownian case (see Aldous \cite{Al1, Al2}). Stable trees (and more generally L\'evy trees) are the scaling limit of Galton-Watson trees (see \cite{DuLG} Chapter 2 and \cite{Du2}). Various geometric and distributional properties of L\'evy trees (and of stable trees, consequently) have been studied in \cite{DuLG2} and in Weill \cite{Weill}.
Stable trees have been also studied in connection with fragmentation processes: see Miermont \cite{Mier03, Mier05}, Haas and Miermont \cite{HaaMi}, Goldschmidt and Haas \cite{GolHaa} for the stable cases and see Abraham and Delmas \cite{AbDel08} for related models concerning more general L\'evy trees. To study Brownian motion on stable trees, D. Croydon in \cite{Croy10} got partial results on balls with exceptionnally small mass measure. 

Before stating the results, let us briefly explain the definition of stable trees before stating the main results of the paper. Let us fix the stable index $\gamma \in (1, 2]$ and let $X= (X_t)_{t\geq 0}$ be a spectrally positive $\gamma$-stable L\'evy process that is defined on the probability space $(\Omega, \cF, \bP)$. More precisely, we suppose that $\bE [\exp (-\lambda X_t)] = \exp (t\lambda^\gamma)$, $t, \lambda \in [0, \infty)$. Note that $X$ is a Brownian motion when $\gamma= 2$ and we shall refer to this case as to the {\it Brownian case}. As shown by Le Gall and Le Jan \cite{LGLJ1} (see also \cite{DuLG} Chapter 1), there exists a {\it continuous process} $H= (H_t )_{ t \geq 0}$ such that for any $t \in [0, \infty)$, the following limit holds true in probability
\begin{equation}
\label{Hlimit}
H_t:=\lim_{\varepsilon\to 0} \frac{1}{\varep}\int_0^t {\bf 1}_{\{I^s_t<X_s<I^s_t+\varep\}}\,ds .
\end{equation}
Here $I^s_t$ stands for $\inf_{s\leq r\leq t} X_r$.
The process $H$ is the $\gamma$-{\it stable height process}. Note that in the Brownian case, $H$ is simply a reflected Brownian motion. Theorems 2.3.2 and 2.4.1 in \cite{DuLG} show that $H$ is the scaling limit of the contour function (or the depth-first exploration process) of a i.i.d.$\;$sequence of Galton-Watson trees whose offspring distribution is in the domain of attraction of a $\gamma$-stable law.

  As in the discrete setting, the process $H$ encodes a family of continuous trees: each excursion of $H$ above $0$ corresponds to the exploration process of a single continuous tree of the family.
Let us make this statement more precise thanks to excursion theory. Recall that $X$ has unbounded variation sample paths. We set $I_t= \inf_{s\in [0, t]} X_s$, that is the infinimum process of $X$. Basic results on fluctuation theory (see Bertoin \cite{Be} VII.1) entail that $X-I$ is a strong Markov process in $[0, \infty)$ and that $0$ is regular for
$(0, \infty)$ and recurrent with respect to this Markov process. Moreover, $-I$ is a local time at $0$ for $X-I$ (see Bertoin \cite{Be} Theorem VII.1). We denote by $\baN$ the corresponding {\it excursion
measure of $X-I$ above $0$}. We denote by $(a_j, b_j)$, $j\in  \cI$, the excursion intervals of $X-I$ above $0$, and by $X^j = X_{(a_j + \cdot )\wedge b_j}-I_{a_j}$, $j\in \cI$, the corresponding excursions.
Next, observe that if $t \in (a_j, b_j)$, the value of $H_t$ only depends on $X^j$. Moreover, one can show that $\bigcup_{^{j\in \cI}} (a_j, b_j)= \{ t \geq 0: H_t >0 \}$.
This allows to define the height process under $\baN$ as a certain measurable function $H(X)$ of $X$, that we simply denote by $(H_t)_{ t\geq 0}$.  For any $j\in  \cI$, we then set $H^j =H_{(a_j+\cdot )\wedge b_j} $ and the point measure
\begin{equation}
\label{Poissheight}
\sum_{j\in \cI} \delta_{(-I_{a_j},H^j)}
\end{equation}
is distributed as a Poisson point measure on $[0, \infty)\times C([0, \infty), \bR)$ with intensity $\ell  \otimes \baN$, where $\ell$ stands for the Lebesgue measure on $[0, \infty)$. Note that $X$ and $H$ under $\baN$ are paths with the same lifetime given by
$$\zeta : = \inf \{ t \in [0, \infty) :  \forall s \in [t, \infty) \, , \, H_s = H_t  \} \;.$$
Standard results in fluctuation theory imply that $0<\zeta< \infty$, $\baN$-a.e.$\;$and that
$$ \baN \big( 1-e^{-\lambda \zeta} \big)= \lambda^{1/\gamma} \; , \; \lambda\in [0, \infty)\; .$$
Thus, $ \baN (\zeta \in dr)= C\,  r^{-\frac{1}{\gamma} -1}
\ell (dr)$,
where $1/C= \gamma \Gamma (1-\frac{1}{\gamma})$ (here, $\Gamma $ stands for Euler's Gamma function).
Note that $\baN$-a.e.$\, H_t >0$ iff $t \in (0, \zeta)$ and $H_0=H_t=0$, for any $t \in [\zeta , \infty)$.
We refer to Chapter 1 in \cite{DuLG} for more details.

 The excursion $(H_t)_{0\leq t \leq \zeta}$ under $\baN$ is the depth-first exploration process of a continuous tree that is defined as the following metric space: for any $s,t \in [0, \zeta]$, we set
\begin{equation}
\label{treedist}
b(s,t)= \min_{s\wedge t \leq r \leq s \vee t} H_r \quad {\rm and} \quad
d(s, t)= H_t+H_s -2 b(s,t) \; .
\end{equation}
The quantity $b(s,t)$ is the height of the branching point between the vertices visited at times $s$ and $t$ and $d(s,t)$ is therefore the distance in the tree of these vertices.
We easily show that $d$ is a pseudo-metric and we introduce the equivalence relation $\sim$ on $[0, \zeta]$ by setting $s\sim t $ iff $d(s, t)= 0$. We then define the
{\it $\gamma$-stable tree} as the quotient metric space
$$ \big( \cT , d \big) = \big( [0, \zeta]/ \! \! \sim \, , \,  d \, \big)  \; .$$
We denote by ${\rm p}: [0, \zeta ] \rightarrow \cT$ the {\it canonical projection}. It is easy to see that ${\rm p}$ is continuous. Thus the $\gamma$-stable tree $(\cT, d)$ is a (random) connected compact metric space.
More precisely, Theorem 2.1 in \cite{DuLG2} asserts that $(\cT, d)$ is a {\it $\bR$-tree}, namely a metric space such that the following holds true for any $\sigma , \sigma^\prime \in \cT$.
\begin{description}
\item[{\bf (a)}] There is a unique isometry
$f_{\sigma,\sigma^\prime}$ from $[0,d(\sigma ,\sigma^\prime)]$ into $\cT$ such
that $f_{\sigma ,\sigma^\prime }(0)=\sigma $ and $f_{\sigma ,\sigma^\prime }(
d (\sigma ,\sigma^\prime))=\sigma^\prime$. We set $\lgeo \sigma,\sigma^\prime \rgeo =f_{\sigma,\sigma^\prime }([0, d (\sigma ,\sigma^\prime)])$ that is the geodesic joining $\sigma $ to $\sigma^\prime$.
\item[{\bf (b)}]   If $g: [0, 1] \rightarrow \cT$ is continuous injective, then $g([0,1])=\lgeo g(0) , g(1) \rgeo$.
\end{description}
We refer to Evans \cite{EvStF} or to Dress, Moulton and Terhalle \cite{DMT96} for a detailed account on $\bR$-trees.
An intrinsic approach of continuous trees has been developped by Evans Pitman and Winter in \cite{EvPitWin} (see also \cite{DuLG2}): Theorem 1 in \cite{EvPitWin} asserts that the set of isometry classes of compact $\bR$-trees equipped with the Gromov-Hausdorff distance (see Gromov \cite{Gro}) is a Polish space. This geometric point of view has been used in \cite{DuLG2} and \cite{DuWi1} to study L\'evy trees.

 We need to introduce two additional features of $\gamma$-stable trees. First, we distinguish a special point $\rho := {\rm p} (0)$ in $\cT$, that is called the {\it root}.
We also equip $\cT$ with the measure $\bm$ that is induced by the Lebesgue measure $\ell $ on $[0, \zeta]$
via the canonical projection ${\rm p}$. Namely, for any Borel subset $A$ of $(\cT, d)$,
$$ \bm (A) = \ell \big( {\rm p}^{-1} (A) \big) \; .$$
The measure $\bm$ is called the {\it mass measure} of the $\gamma$-stable tree $\cT$. Note that
$\bm (\cT)= \zeta$. One can prove that $\bm$ is diffuse and that its
topological support is $\cT$. Moreover,
$\bm$ is carried by the set of {\it leaves of $\cT$} that is the set of the points $\sigma $ such that $\cT \backslash \{ \sigma \}$ remains connected (see \cite{DuLG2} for more details). The measured tree
$(\cT, d , \bm)$ is thus a {\it continuum tree}, as defined by Aldous in \cite{Al2}.

Let us discuss briefly the scaling property of $\cT$. From the scaling property of $X$ and from (\ref{Hlimit}), we see that for any $r \in (0, \infty)$, under $\bP$,
$(r^{\frac{_{\gamma-1}}{^{\gamma}}}H_{t/r})_{t\geq 0}$ has the same distribution as $H$. Then, by (\ref{Poissheight}), $(r^{\frac{\gamma-1}{\gamma}}H_{t/r})_{0\leq t\leq r\zeta}$ under $r^{-\frac{1}{\gamma}} \baN$ has the same "distribution" as $(H_t)_{0\leq t\leq \zeta}$ under $\baN$.
Thus, the rescaled measured tree
$(\cT \, , \,  r^{\frac{\gamma-1}{\gamma}} \! d \, , \, r  \bm)$ under $r^{-\frac{1}{\gamma}} \baN$ has the same "distribution" as $(\cT, d, \bm)$ under $\baN$.
This allows to define for any
$r \in (0, \infty)$, a probability distribution $\baN ( \, \cdot \, | \, \zeta = r)$ on $C([0, \infty), \bR)$ such that
$r \mapsto \baN( \, \cdot \, | \, \zeta = r)$ is weakly continuous and such that
$$ \baN = \int_0^\infty \!\!\! \baN( \, \cdot \, | \, \zeta = r)\,  \baN( \zeta \in dr) \; .$$
Moreover, $(r^{\frac{\gamma-1}{\gamma}} H_{t/r})_{0\leq t \leq r} $ under $\baN  (\, \cdot \, | \, \zeta = 1) $ has the same distribution as $(H_{t})_{0\leq t\leq r}$
under $\baN( \, \cdot \, | \, \zeta = r)$. Since $\bm (\cT)= \zeta$, the tree $\cT$ under $\baN  (\, \cdot \, | \, \zeta = 1) $ is interpreted as the $\gamma$-stable tree conditioned to have total mass equal to $1$ and it is simply called the {\it normalised $\gamma$-stable tree}. When $\gamma= 2$, it corresponds (up to a scaling constant) to Aldous Continuum Random Tree as defined in \cite{Al1} (see also Le Gall
\cite{LG2} for a definition via the normalised Brownian excursion). The normalised $\gamma$-stable tree is the weak limit when $n$ goes to infinity of a rescaled Galton-Watson trees conditionned
to have total size $n$ and whose offspring distribution belongs to the domain of attraction of a $\gamma$-stable law: see Aldous \cite{Al2} for the Brownian case and see \cite{Du2} for the general case.

  The mass measure $\bm$ is in some sense the most spread out measure on $\cT$ and it plays a crucial role in the study of stable trees.
For instance Theorem 1.1 in \cite{Du10} asserts that for any $\gamma \in (1, 2]$, $\baN$-a.e.$\;$the mass measure $\bm$ is equal to a deterministic constant times the $g_\gamma$-packing measure where the gauge function is given by
\begin{equation}
\label{Gauge}
   g_\gamma (r)= \frac{r^{\frac{\gamma}{\gamma -1}} }{(\log \log 1/r)^{\frac{1}{\gamma-1}} }\, , \quad r \in (0, e^{-1} ) \; .
\end{equation}
Actually, this result holds true for general L\'evy trees (with a more involved gauge function).
Here, the power exponent $\gamma/ (\gamma-1)$ reflects the scaling property. This value
is also equal to the packing dimension of $\cT$, and to its Hausdorff and its box counting dimensions (see \cite{DuLG2}).
The function $g_\gamma$ is also the {\it lower density of $\bm$} at typical points. More precisely, denote by $B (\sigma, r)$ the open ball in $\cT$ with center $\sigma \in \cT$ and radius $r \in (0, \infty)$. Then, Theorem 1.2 in \cite{Du10}  asserts that there exists a constant $C_\gamma \in (0, \infty)$ such that
\begin{equation}
\label{lowdens}
\textrm{$\baN$-a.e.$\;$for $\bm$-almost all $\sigma$,} \quad \liminf_{r\rightarrow 0} \frac{\bm \big( B(\sigma, r) \big)}{g_\gamma (r)} = C_\gamma \; .
\end{equation}
Theorem 1.2 in \cite{Du10} also holds true 
for general L\'evy trees and the constant is unknown. However,
in the stable cases, we are able to compute explicitely $C_\gamma$, as shown by  the following proposition.
\begin{proposition}
\label{Constant} For any $\gamma \in (1, 2]$, $C_\gamma = \gamma-1$.
\end{proposition}
We also discuss the balls with exceptionally small mass measure. More precisely, we investigate the  behaviour of $ \inf_{\sigma \in \cT} \bm \big( B(\sigma, r) \big)$ when $r$ goes to $0$.
Our main result is the following.
\begin{theorem}
\label{thin} For any $\gamma \in (1, 2]$, we set
\begin{equation}
\label{smallgauge}
 f_\gamma (r)= \frac{r^{\frac{\gamma}{\gamma -1}} }{
( \log 1/r)^{\frac{1}{\gamma-1}} } \; , \quad r \in (0, 1) \; .
\end{equation}
Then, there exist $k_\gamma , K_\gamma \in (0, \infty) $ such that $\baN$-a.e.
 \begin{equation}
\label{thindens}
k_\gamma \;  \leq \;
\liminf_{r\rightarrow 0} \frac{1}{f_\gamma (r)} \inf_{\sigma \in \cT} \bm \big( B(\sigma, r) \big) \;   \leq \;  \limsup_{r\rightarrow 0} \frac{1}{f_\gamma (r)} \inf_{\sigma \in \cT} \bm \big( B(\sigma, r) \big) \; \leq  \; K_\gamma \; .
\end{equation}
\end{theorem}
To study Brownian motion on stable trees, D. Croydon in 
Proposition 5.1 \cite{Croy10} states a partial lower bound for $\inf_{\sigma \in \cT} \bm \big( B(\sigma, r) \big) $ that is sufficient to his purpose (but that does not provide the right scale function).

  When $1< \gamma <2$, there is no exact upper density of $\bm$ at typical points (see Proposition 1.9 in \cite{Du9}). Moreover, Theorem 1.10 in \cite{Du9} shows that $\cT$ has no exact Hausdorff measure whose gauge function is regularly varying.
We also strongly believe that there is no exact asymptotic function for
$r \mapsto  \sup_{\sigma \in \cT} \bm \big( B(\sigma , r)\big)$, when $r$ goes to $0$, but we will not consider this problem in this paper.

In the Brownian case, Theorem 1.1 in \cite{DuLG3} asserts that $\baN$-a.e.$\;$the mass measure $\bm$ is equal to a deterministic constant times the $g$-Hausdorff measure where the gauge function $g$ is given by
$$g (r)= r^{2}\log \log 1/r \; , \quad r \in (0, e^{-1} ) \; .$$
In the proof of Theorem 1.1 \cite{DuLG3}, it is proved that $g$ is the upper density of $\bm$ at typical points. In this paper we obtain the following specific constant.
\begin{proposition}
\label{constant} Consider the Brownian case: $\gamma= 2$. Then
\begin{equation}
\label{uppdens}
\textrm{$\baN$-a.e.$\;$for $\bm$-almost all $\sigma$,} \quad \limsup_{r\rightarrow 0} \frac{\bm \big( B(\sigma, r) \big)}{g (r)} = \frac{4}{\pi^2} \; .
\end{equation}
\end{proposition}
Moreover, in the Brownian case, the balls with exceptionally large mass have also an exact asymptotic function as shown by the following theorem.
\begin{theorem}
\label{ThickB} Consider the Brownian case: $\gamma= 2$. Let us set
$$f(r)= r^2 \log 1/r \; , \quad r\in (0, e^{-1}) \; .$$
Then, there exist $k , K \in (0, \infty) $ such that $\baN$-a.e.
 \begin{equation}
\label{ThicBdens}
k  \; \leq \;
\liminf_{r\rightarrow 0} \frac{1}{f (r)} \sup_{\sigma \in \cT} \bm \big( B(\sigma, r) \big)  \; \leq \;
 \limsup_{r\rightarrow 0} \frac{1}{f (r)} \sup_{\sigma \in \cT} \bm \big( B(\sigma, r) \big)\;  \leq \; K \; .
\end{equation}
\end{theorem}
Observe that, in the Brownian case, Theorem \ref{thin} and Theorem \ref{ThickB} immediately imply the following result.
\begin{corollary}
\label{coincB} Consider the Brownian case: $\gamma= 2$. Then, there are two constants $c, C\in (0, \infty)$ such that 
$$\textrm{$\baN$-a.e.} \quad  \exists r_0  \!\in \!(0, \infty) \; : \quad \forall r \! \in \! (0, r_0)  , \, \forall \sigma \! \in \! \cT  , \quad  \frac{c}{\log 1/r}  \leq   r^{-2} \bm \big( B(\sigma , r)\big) \leq    C \log 1/r  . $$
\end{corollary}

  Note that Proposition \ref{Constant}, Theorem \ref{thin}, Proposition \ref{constant}, Theorem \ref{ThickB} and Corollary \ref{coincB} hold true under the normalised law $\baN( \, \cdot \, | \, \zeta= 1)$.

\medskip

 The paper is organised as follows. Section \ref{basicsec} recalls useful technical results on the height process and basic geometric properties of stable trees. Section \ref{tailsec} is devoted to the tail estimates of the mass measure of specific subsets of stable trees.
 Section \ref{proofsec} is devoted to the proof of the results.

\section{Preliminaries and basic results.}
\label{basicsec}
\subsection{Results on the stable height process.}
\paragraph{Local times of the height process.}
Let $H$ be the $\gamma$-stable height process under its excursion measure $\baN$ as defined in the introduction. It is possible to define the local times of $H$ under the excursion measure $\baN$
as follows. For any $b \in (0, \infty)$, let us  set $v(b)=\baN ( \sup_{^{t \in [0, \zeta ]}} H_t > b )$. The continuity of $H$ under $\bP$ and the Poisson decomposition (\ref{Poissheight}) obviously imply that $v(b) < \infty$, for any $b >0$.
 It is moreover clear that $v$ is non-increasing and $\lim_{b \rightarrow \infty} v (b)= 0$.
 For every $a\in (0, \infty)$, we then define a continuous increasing process $(L^a_t)_{0\leq t \leq  \zeta}$, such that for every
$b \in (0,\infty)$ and for any $t\in [0, \infty)$, one has
\begin{equation}
\label{localapprox}
\lim_{\varepsilon \rightarrow 0} \,
\baN \left(  \un_{\{\sup H>b \}}\;\sup_{ 0\leq s \leq t\wedge \zeta} \left| \frac{1}{\varepsilon} \int_0^s
dr
\un_{\{ a-\varepsilon< H_r \leq a\}} -L_s^a \right| \right) =0.
\end{equation}
See \cite{DuLG} Section 1.3 for more details. The process $(L^a_t)_{0\leq t \leq  \zeta}$ is the {\it $a$-local time of the height process}. For any $a, \lambda, \mu \in [0, \infty)$, we set
\begin{equation}
\label{kappaexc}
 \kappa_a (\lambda, \mu) := \baN \left( 1- e^{ -\mu L^a_\zeta - \lambda \int_0^a \un_{\{ H_t < a \}} dt} \right).
\end{equation}
The function $\kappa$
is the Laplace exponent of a specific additive functional of a $\gamma$-stable Continuous States Branching Process ($\gamma$-stable CSBP, for short).
An elementary result on CSBPs, whose proof can be found in Le Gall \cite{LG99} Section II.3,
entails that
$a\mapsto \kappa_a (\lambda , \mu)$ is the unique solution of the following ordinary
differential equation:
$$\frac{\partial \kappa_a}{\partial a}  (\lambda ,\mu )= \lambda -\kappa_a (\lambda, \mu)^\gamma \quad {\rm  and}  \quad \kappa_0 (\lambda, \mu) = \mu \; .$$
For more details, see \cite{DuLG3} Section 4 page 405 or \cite{Du9} Section 2.3 page 106. We note the following: if $\mu =\lambda^{1/\gamma}$, then $\kappa_a (\lambda , \mu)= \lambda^{1/\gamma}$. If $\mu <\lambda^{1/\gamma}$, (resp.$\;$if $\mu >\lambda^{1/\gamma}$), then $a\mapsto \kappa_a (\lambda , \mu)$ is increasing (resp.$\;$decreasing). A simple change of variable implies that $\kappa$ satisfies the following integral equation
\begin{equation}
\label{equakappa}
\int_{\mu}^{\kappa_a (\lambda, \mu) } \frac{du}{\lambda-u^\gamma} = a , \qquad \textrm{for any $a, \lambda , \mu \in [0, \infty )$ such that $\mu \neq \lambda^{\frac{1}{\gamma}}$. }
\end{equation}
As a consequence, we get
\begin{equation}
\label{kappabranch}
\kappa_{a+b} (\lambda, \mu)= \kappa_a \big( \lambda, \kappa_b (\lambda, \mu ) \big) \, , \quad a, b, \lambda , \mu \in [0, \infty) \; .
\end{equation}
We also derive from (\ref{equakappa}) the following scaling property:
\begin{equation}
\label{kappascale}
c^{\frac{1}{\gamma-1}} \kappa_a \big( c^{-\frac{\gamma}{\gamma-1}} \lambda ,  c^{-\frac{1}{\gamma-1}}    \mu \big) = \kappa_{a/c} (\lambda, \mu) \; , \quad a,  c , \lambda, \mu \in [0, \infty) \; .
\end{equation}
When $\gamma \neq 2$ and $\lambda >0$, it seems difficult to compute $\kappa$ explicitly. However, when $\lambda= 0$, (\ref{equakappa}) implies that for any $a, \mu \in [0, \infty)$, 
\begin{equation}
\label{explisolbranch}
\kappa_a (0, \mu) =\baN \big( 1- e^{ -\mu L^a_\zeta } \big)=  \Big( (\gamma\! -\!1) a  +\frac{1}{{\mu^{\gamma -1}}}  \Big)^{-\frac{1}{\gamma-1}}\, .
\end{equation}

  It is convenient to interpret these quantities in terms of the $\gamma$-tree $\cT$.
For any $a \in (0, \infty)$, first define the {\it $a$-local time measure} $\ell^a$ as the measure induced by $dL^a_{\cdot}$ via the canonical projection ${\rm p}: [0, \zeta] \rightarrow \cT$. Namely,
  $$ \langle \ell^a , f \rangle = \int_0^\zeta f({\rm p} (s)) \, dL^a_s \; , $$
for any positive measurable application $f$ on $\cT$. Here $dL^a_{\cdot}$ stands for the Stieltjes measure associated with the non-decreasing function $s \mapsto L^a_s $.
Note that the topological support of $\ell^a$ is included in  the $a$-level set
\begin{equation}
\label{aleveldef}
\cT (a)= \{ \sigma \in \cT: d(\rho ,\sigma)= a \} \; .
\end{equation}
Next, we set
\begin{equation}
\label{heightdef}
\Gamma (\cT)= \sup_{\sigma \in \cT} d( \rho , \sigma)
\end{equation}
that is the {\it total height of $\cT$}. Then, observe that
\begin{equation}
\label{identif}
\langle \ell^a\rangle= L^a_\zeta \, , \quad \bm \big( B(\rho , a) \big)= \int_0^\zeta \!\! \un_{\{ H_t <a \}} dt \quad \textrm{and} \quad \Gamma (\cT) = \sup_{t\in [0, \zeta]} H_t \; ,
\end{equation}
where $\langle \ell^a\rangle$ stands for the total mass of $\ell^a$. This implies
\begin{equation}
\label{kappatree}
\kappa_a (\lambda, \mu) = \baN \Big( 1- e^{ -\mu \,  \langle \ell^a\rangle - \lambda\,  \bm ( B(\rho , a) )} \Big).
\end{equation}
We recall from \cite{DuLG} Chapter 1 (proof of Corollary 1.4.2 page 41) that 
\begin{equation}
\label{alevelNae}
\textrm{ $\baN$-a.e.} \quad  \un_{\{ \sup H >a \}}= \un_{\{ L^a_\zeta \neq 0\}} \; .
\end{equation}
By letting $\mu$ go to $\infty$ in (\ref{explisolbranch}), we get
\begin{equation}
\label{heightform}
v(a)= \baN (\Gamma (\cT) >a ) = \baN ( \langle \ell^a\rangle\neq 0) = \big( \, (\gamma \!- \!1)a \big)^{-\frac{1}{\gamma -1}} \; .
\end{equation}
For any $a \in (0, \infty)$, we next define the probability measure $\baN_a$ by setting
\begin{equation}
\label{Nadef}
\baN_a= \baN \big( \, \cdot \, \big| \,  \Gamma (\cT) >a    \big)= \baN \big( \, \cdot \, \big| \, \langle \ell^a \rangle  \neq 0\big) \; .
\end{equation}
This combined with (\ref{explisolbranch}) implies that
\begin{equation}
\label{Nabranch}
\baN_a \big(  e^{-\mu \langle \ell^a \rangle }\big)= 1- \Big( 1+ \frac{1}{(\gamma\!- \!1)a\mu^{\gamma-1}} \Big)^{-\frac{1}{\gamma-1}} \; .
\end{equation}
By differentiating this equality at $\mu= 0$, one gets
\begin{equation}
\label{Nabranchmean}
\baN_a \big(  \langle \ell^a \rangle \big)=\big( \, (\gamma \!- \!1)a \big)^{\frac{1}{\gamma -1}} =  \frac{1}{v(a)} \; .
\end{equation}

\paragraph{The branching property.}
We now  describe the distribution of excursions of the height process above a given level (or equivalently of the corresponding stable subtrees above this level). We fix $a\in (0, \infty)$ and we denote by $(l^{_a}_{^j}, r^{_a}_{^j})$, $j\in \cJ_a$,
the connected components of the open set $\{t \in (0, \zeta) :H_t>a\}$. For any $j\in \cJ_a$,
denote by $H^{_{a, j}}_{^{.}}$ the corresponding excursion of $H$ that is defined by $H^{_{a, j}}_s=H_{(l^a_j+s)\wedge r^a_j}-a$, $ s\in [ 0, \infty) $.

  This decomposition is interpreted in terms of the tree as follows. Denote the closed ball in $\cT$ with center $\rho$ and radius $a$ by $\bar{B} (\rho , a)$. Observe that the connected components of the open set $\cT \backslash \bar{B} (\rho , a)$ are the subtrees $\cT_{^j}^{_{a, o}}:= {\rm p}((l^{_a}_{^j}, r^{_a}_{^j}))$, $j \in \cJ_a$. The closure $\cT_{^j}^{_a}$ of $\cT_{^j}^{_{a, o}}$ is simply $\{ \sigma^{_a}_{^j} \} \cup \cT_{^j}^{_{a, o}}$, where $\sigma^{_a}_{^j} = {\rm p}(l^{_a}_{^j})= {\rm p}(r^{_a}_{^j})$, that is the point in the $a$-level set $\cT(a)$ at which $\cT^{_{a, o}}_{^j}$ is grafted.
Observe that the rooted measured tree $\big( \, \cT^{_a}_{^j}, d, \sigma^{_a}_{^j} , \bm (\, \cdot \, \cap \cT^{_a}_{^j}) \, \big)$ is isometric to the tree coded by $H^{_{^{a, j}}}_{^{.} }$.

We define the following point measure on $[0, \infty) \times C([0, \infty) , \bR)$:
\begin{equation}
\label{branchprop}
\cM_a (dxdH)=  \sum_{^{j\in \cJ_a}} \delta_{( L^a_{l^a_j} ,H^{a,j}  )}
\end{equation}
For any $s \in [0, \infty)$, we also set $\tilde{H}^a_s=H_{\tau^a_s}$, where the time-change $\tau^a_s$ is given by
$$ \tau^a_s=\inf  \big\{  t\geq 0:\int_0^t dr\,\un_{\{H_r\leq a\}}>s \big\}\; , \quad s\in (0, \infty)\; .$$
The process $\tilde{H}^a=(\tilde{H}^a_s)_{ s \geq 0}$ is the {\it height process below $a$} and the rooted compact $\bR$-tree
$(\bar{B} (\rho , a), d, \rho)$ is isometric to the tree coded by $\tilde{H}^a$. Let $\cG_a$ be the
sigma-field generated
by $\tilde{H}^a $ augmented by the $\baN$-negligible sets. From the approximation
(\ref{localapprox}), it follows that $L^{_a}_{^\zeta}$ is measurable with respect
to $\cG_a$. Recall notation $\baN_a$ from (\ref{Nadef}).

\medskip

\noi
{\it The branching property at level $a$} then asserts that under $\baN_a$, conditionally given $\cG_a$, $\cM_a $ is distributed as a Poisson point measure with intensity
$\un_{[0,L^{a}_{\zeta} ]}(x) \ell (dx) \otimes  \baN (dH)$. 

\medskip

\noi
For more details, we refer to Proposition 1.3.1 in \cite{DuLG} or to the proof of Proposition 4.2.3 \cite{DuLG} (see also Theorem 4.2 \cite{DuLG2}). We apply the branching property to prove the following lemma.
\begin{lemma}
\label{kappatriangle} For any $a, b, \lambda \in (0, \infty)$, we set
$$\Phi_{a, b} (\lambda)= \baN_a \big( e^{ -\lambda \bm (B(\rho , a+b)) } \big) \; .$$
Then,
$$\Phi_{a, b} (\lambda) = v(a)^{-1} \big( \kappa_a  (\lambda, \infty ) -\kappa_{a+b} ( \lambda, 0) \big) \; , $$
where $\kappa_a (\lambda, \infty) $ stands for $\lim_{\mu \rightarrow \infty} \kappa_a (\lambda, \mu) $, which is well-defined and finite.
\end{lemma}
\noi
{\bf Proof:} first note that
\begin{eqnarray*}
\bm (B(\rho , a+b)) &= & \int_0^\zeta \!\! \un_{\{ H_r < a+b \}} dr = \int_0^\zeta \!\! \un_{\{ H_r < a \}} dr+\int_0^\zeta \!\! \un_{\{ a\leq H_r < a+b \}} dr \\
&= &  \int_0^\zeta \!\! \un_{\{ H_r < a \}} dr + \sum_{^{j\in \cJ_a}} \int_0^{\zeta^a_j} \!\! \un_{\{ H^{a, j}_r < b \}} dr \;  ,
\end{eqnarray*}
where $\zeta^a_j= r^a_j-l^a_j$ stands for the lifetime of $H^{a, j}$. The branching property then implies that
\begin{eqnarray*}
\baN_a \big( e^{ -\lambda \bm (B(\rho , a+b)) } \, \big| \, \cG_a  \big) &=& e^{ -\lambda \int_0^\zeta \!\! \un_{\{ H_r < a \}} dr } \;
\exp\! \big( \!\! - \! L^a_\zeta \baN \big(1-e^{-\lambda \int_0^\zeta \!\! \un_{\{ H_r < b \}} dr} \big)     \big) \\
& =& e^{- \kappa_b (\lambda , 0) L^a_\zeta -\lambda \int_0^\zeta \!\! \un_{\{ H_r < a \}} dr } \; .
\end{eqnarray*}
Monotone convergence implies
\begin{eqnarray*}
 \baN_a \big(  e^{- \kappa_b (\lambda , 0) L^a_\zeta -\lambda \int_0^\zeta \!\! \un_{\{ H_r < a \}} dr } \big)
 \!\!\!\! & = & \!\!\!\! \lim_{\mu \rightarrow \infty} v(a)^{-1} \baN \big(  (1 -e^{-\mu L^a_\zeta}  )  e^{- \kappa_b (\lambda , 0) L^a_\zeta -\lambda \int_0^\zeta \!\! \un_{\{ H_r < a \}} dr } \big) \\
  \!\!\!\! & = & \!\!\!\! \lim_{\mu \rightarrow \infty} v(a)^{-1} \big( \, \kappa_{a} \big( \lambda , \mu \!+\! \kappa_b (\lambda , 0) \big) - \kappa_a (\lambda , \kappa_b (\lambda, 0) \big)\, \big) \; ,
\end{eqnarray*}
which easily implies the desired result thanks to (\ref{kappabranch}). \cqfd

\medskip

From Lemma \ref{kappatriangle} and the scaling property (\ref{kappascale}), we get
$$ \Phi_{a, b} (\lambda) = \Phi_{1, b/a} \big( a^{\frac{\gamma}{\gamma-1}} \lambda\big) \; , $$
which implies that for any $a, b\in (0, \infty) $,
\begin{equation}
\label{trianscale}
\textrm{$\bm \big( B (\rho , a+ b) \big)$ under $\baN_a$} \; \overset{\textrm{(law)}}{=} \;   \textrm{$a^{\frac{\gamma}{\gamma-1}} \bm \big( B (\rho , 1+ \frac{_a}{^b}) \big)$ under $\baN_1$} \; .
\end{equation}

\paragraph{Spinal decomposition.}
 We recall another decomposition of the height process that is proved in \cite{DuLG} Chapter 2 
(see \cite{DuLG2} for a more specific statement  and see \cite{DuLG4} for further applications). This decomposition is used in the proof of Proposition \ref{Constant} and Proposition \ref{constant}.

 Let $h : [0, \infty) \rightarrow [0, \infty)$ be a continuous function with compact support. Let us assume for clarity that $h(0) >0$. We view $h$ as the depth-first exploration process of a tree. Thus, the exploration starts at a vertex with height $h(0) >0$ that we call the initial vertex. We obtain the subtrees grafted along the ancestral line of the initial vertex as follows: set $\underline{h} (s)= \inf_{[0, s]} h$ and denote by $(l_i, r_i)$, $i\in \cI(h)$, the excursion intervals of $h-\underline{h} $ away from $0$ that are the connected components of the open set $\{ s \in [0, \infty): h(s)-\underline{h} (s) >0 \}$. 
For any $i \in \cI (h)$, set
$$h^i (s)  = \big(\, (h- \underline{h}) ( (l_i +s) \wedge r_i) \, \big)_{  s \geq 0} \; .$$
Then, the subtrees along the ancestral line of the initial vertex are coded by the functions $(h^i\, ; \, i\in \cI(h))$, and the tree coded by $h^i$ is grafted at distance $h(0)-h(l_i)$ from the initial vertex.
We next define $\cN(h)$ as the point measure on $[0, \infty) \times C ([0, \infty) , \bR)$ given by
$$ \cN (h)= \sum_{ i\in \cI (h)} \delta_{( h(0)- h(l_i) \, , \,  h^i ) } \; .$$

   Recall that $H$ stands for the $\gamma$-height process under its excursion measure $\baN$. For any $t \in (0, \zeta)$, we set $\hat{H}^{t}:= ( H_{(t- s)_+})_{ s \geq 0} $; here, $(\, \cdot )_+$ stands for the positive part function. We also set $\check{H}^{t}:= (H_{(t+s)\wedge \zeta} )_{ s \geq 0}$, and we define the random point measure $\cN_t$ on  $[0, \infty) \times C ([0, \infty) , \bR)$ by
\begin{equation}
\label{spinaldef}
\cN_t =\cN( \hat{H}^{t})+ \cN ( \check{H}^{t}) := \sum_{j \in \cJ_t } \delta_{( r^t_j \, , \,  H^{(t),j} )} \; .
\end{equation}
This point measure records the subtrees grafted along the ancestral line of the vertex visited at time $t$ in
the coding of $\cT$ by $H$. Namely, set $\sigma = {\rm p}(t) \in \cT$. Then, the
geodesic $\lgeo \rho, \sigma \rgeo$ is the ancestral line of $\sigma$. Denote by
$\cT_{^j}^{_o}$, $j \in \cJ$, the connected components of the open set $\cT   \backslash \lgeo \rho , \sigma \rgeo$ and denote by $\cT_{{j}}$ the closure of
$\cT_{^j}^{_o}$. Then, there exists a point $\sigma_j \in \lgeo \rho , \sigma \rgeo$ such that $\cT_{{j}} = \{ \sigma_{{j}} \} \cup \cT_{^j}^{_o}$. 
The specific coding of $\cT$ by $H$ entails that for any $j \in \cJ$ there exists a unique $j^\prime \in \cJ_{t}$ such that $ d(\sigma ,\sigma_{{j}})= r^{t}_{^{j^\prime}} $ and such that the rooted compact $\bR$-tree $(\cT_{{j}}, d, \sigma_{{j}})$ is isometric to the tree coded by $H^{_{^{(t),j^\prime}}}_{^{.}}$.

 The law of $\cN_t$ when $t$ is chosen "at random" according to the Lebesgue measure is given as follows. To simplify the argument, we assume that the random variables we mention are defined on the same probability space $(\Omega, \cF, \bP)$. Let $(U_t)_{ t \geq 0}$ be a ($\gamma-1$)-stable subordinator
with initial value $U_0= 0$, and whose Laplace exponent is given by $-\log\bE [ \exp (-\lambda U_1)] = \gamma \lambda^{\gamma-1}$. Let
\begin{equation}
\label{Nstardef}
\cN^* = \sum_{j \in \cI^*} \delta_{ (  r^*_j , \, H^{*j} )  }
\end{equation}
be a random point measure on $[0, \infty) \times C ([0, \infty), \bR)$ such that a regular version of the law of $\cN^*$ conditionally given $U$ is that of a Poisson point measure with intensity $dU_r \otimes \baN (dH)$. Here, $d U_r$ stands for the (random) Stieltjes measure associated with the non-decreasing path $r \mapsto U_r$. For any $ a \in (0, \infty)$, we also set
\begin{equation}
\label{Nstaradef}
\cN^*_a =  \sum_{j \in \cI^*} \un_{[0, a]} (r^*_j) \, \delta_{ (  r^*_j , \, H^{*j} )  } .
\end{equation}
By Lemma 3.4 in \cite{DuLG2}, for any nonnegative measurable functional $F$, 
\begin{equation}
\label{ancdecomp}
\baN \left( \int_{0}^\zeta \!\!\! F \big( H_t,  \cN_t  \big) \, dt \right) =
\int_0^\infty \!\! \bE \left[ F (a, \cN^*_a ) \right]  \, da \; .
\end{equation}
We shall refer to this identity as to the {\it spinal decomposition of $H$ at a random time}.

\medskip

We use the spinal decomposition to compute the law of the mass measure of random balls in $\cT$.
To that end, we first fix $t \in (0, \zeta )$ and we express $\bm \big( B({\rm p}(t), r) \big)$ in terms of  $\cN_t$ as follows. First, recall from (\ref{treedist}) the definition of $b(s,t) $ and $d(s, t)$. Note that if $H_s= b(s,t)$ with $s \neq t$, then ${\rm p}(s)$ cannot be a leaf of $\cT$.
Let us fix a radius $r$ in $[0, H_t]$. Since the leaves of $\cT$ have zero $\bm$-measure, we get
$$ \bm \big( B ({\rm p}(t) , r  ) \big) = \int_{0}^\zeta \un_{\{ d(s,t) <  r \}} ds = \int_0^\zeta
\un_{\{ 0 < H_s  - b(s,t ) < r -  H_t+b(s,t)     \}}  ds . $$
The definition of  $(\cN(\hat{H}^t), \cN(\check{H}^t))$ then entails
\begin{equation}
\label{massballfixt}
\bm \big( B({\rm p} (t) , r  ) \big) = \sum_{j \in \cJ_t} \un_{[ 0, r]  } (r^t_j)\cdot  \int_0^{\zeta^t_j} \un_{\{ H^{(t),j}_s <  r-r^t_j\}} \, ds ,
\end{equation}
where $\zeta^t_j$ stands for the lifetime of the path $H^{*\,  t,j}$. For any $a \in (0, \infty)$ and for any $r \in [0,a]$, we next set
\begin{equation}
\label{Mstardef}
M^*_r (a) = \sum_{j \in \cI^*} \un_{[0, r\wedge a ]} (r^*_j) \cdot \int_{0}^{\zeta^*_j} \un_{\{ H^{*j}_s \leq r-r^*_j  \}} ds  \; ,
 \end{equation}
where $\zeta^*_j$ stands for the lifetime of the path $H^{*j}$. Then, $(M^*_r (a) )_{0\leq r \leq a} $ is a function of $\cN_a^*$. It  is a cadlag non-decreasing process and the spinal decomposition (\ref{ancdecomp}) implies for any bounded measurable functional  $F: \bbD ([0, \infty) ,
\bR)  \rightarrow \bR$, we have
$$ \baN \left( \int_0^\zeta \!\!  F \! \left(  \big( \bm ( B({\rm p} (t) ,  r ) ) \, \big)_{r \geq 0}  \right) \, dt \right)
= \int_0^\infty \!\! \bE \left[  F \big( (M^*_ r (a) )_{r \geq 0} \big)\right] \, da \, , $$
which can be rephrased as follows in terms of the tree
\begin{equation}
\label{keymass}
 \baN \left( \int_{\cT}  \! F \! \left(  \big( \bm( \!\! \;  B( \sigma , r )) \, \big)_{r \geq 0}  \right) \, \bm (d\sigma) \right)
= \int_0^\infty \!\! \bE \left[  F \big( (M^*_{ r} (a) )_{r \geq 0} \big)\right] \, da \; .
\end{equation}

  For any $0\leq r^\prime\leq r\leq a $, we next set
\begin{equation}
\label{massshelldef}
M^*_{r^\prime, r} (a)=\sum_{j \in \cI^*} \un_{( r^\prime , r  ]} (r^*_j) \int_0^{\zeta^*_j}  \un_{\{ H^{*j}_s
< r- r^*_j \}} \;
\end{equation}
Note that $M^*_r (a)= M^*_{0, r} (a)$. The random variables $M^*_{r^\prime , r} (a)$
play an important role in the proof of Proposition \ref{Constant} and Proposition \ref{constant}.
We gather in the following lemma their basic properties
that are easy consequences of the definition (we refer to Lemma 2.11 and to Remark 2.12 in \cite{Du9} page 115 for more details).
\begin{lemma}
\label{basicMprop} Let us fix $a \in (0, \infty)$. The following holds true.
\begin{itemize}
\item[(i)] Let $(r_n )_{ n \geq 0}$ be a sequence such that $0 < r_{n+1} \leq r_n \leq a$ and $\lim_{n\rightarrow \infty} r_n = 0$. Then, the random variables $(M^*_{r_{n+1}, r_n } (a))_{ n \geq 0}$ are independent.
\item[(ii)] The increments of $r \in [0, a] \mapsto M^*_r (a)$ are not independent. However, for any $0 \leq r^\prime \leq r \leq a $, we have
$$ M_r^* (a) -M^*_{r^\prime} (a)= M^*_{r^\prime , r} (a)+ \sum_{^{j \in \cI^*}} \un_{[0, r^\prime  ]} (r^*_j) \int_0^{\zeta^*_j}  \un_{\{  r^\prime-r^*_j   \leq H^{*j}_s <  r-r^*_j \}} \, .   $$
It implies that $ M_r^* (a)  \geq M_r^* (a)-M_{r^\prime}^* (a)  \geq M^*_{r^\prime , r} (a)$. Note
that $M^*_{r^\prime , r} (a)$ is independent of $M^*_{r^\prime} (a)$.
\end{itemize}
\end{lemma}
\noi
The law of $M_{r^\prime , r}^* (a)$ is characterised by its Laplace transform:
\begin{equation}
\label{QLapcour}
\bE \big[ \exp (-\lambda  M^*_{r^\prime , r} (a) ) \big]=  1 - \frac{\kappa_{r-r^\prime} (\lambda, 0)^\gamma}{\lambda} \, , \quad \lambda \in [0, \infty) \; .
\end{equation}
\noi
{\bf Proof of (\ref{QLapcour})}: recall that conditionally given $U$, $\cN^*$ is a Poisson point process with intensity $dU_r \otimes \baN$. Thus,
$$ \bE \big[ \exp (-\lambda  M^*_{r^\prime , r} (a) ) \, \big| \, U  \big]= \exp \Big( -\int_{ (r^\prime, r]}
dU_s  \,  \kappa_{r-s } (\lambda, 0) \Big) .$$
We get $\bE [ \exp (-\lambda  M^*_{r^\prime , r} (a) )]=\exp (-\gamma  
\int_{0}^{ r-r^\prime } \kappa_{s}( \lambda, 0)^{\gamma -1}ds  )$.  
By a change of variable based on (\ref{equakappa}), we obtain $\gamma \int_{0}^{ r-r^\prime \, }   \kappa_{s}( \lambda, 0)^{\gamma -1}ds =
\log \lambda - \log \big( \lambda -\kappa_{r-r^\prime} (\lambda , 0)\, \big) $,
which entails (\ref{QLapcour}). \cqfd

\medskip

We also introduce the following notation:
\begin{equation}
\label{Mstarun}
 M_* := M^*_{0, 1} (1)= M^*_1(1)   \; .
\end{equation}
Then, (\ref{QLapcour}) and the scaling property (\ref{kappascale}) imply that for any $0\leq r^\prime \leq r \leq a$,
\begin{equation}
\label{Mstarscale}
(r-r^\prime)^{-\frac{\gamma}{\gamma-1}}   M^*_{r^\prime , r} (a)  \overset{{\rm (law)}}{=}
M_*  \; .
\end{equation}
We see in particular that $r^{-\frac{\gamma}{\gamma-1}}M^*_r (a)$ has the same law as
$M_* $. The tail at $0+$ of the distribution of $M_*$ is studied in Section \ref{tailsec}.

\subsection{Balls and truncated subtrees.}
\label{Balltrianglesec}
Recall that $(H_t)_{0\leq t \leq \zeta }$ stands for the excursion of the $\gamma$-stable process under $\baN$ and that
$(\cT, d)$ is the $\gamma$-stable tree coded by $H$. Recall that for any
$\sigma, \sigma^\prime \in \cT$, $\lgeo \sigma , \sigma^\prime \rgeo$ stands for the (unique)
geodesic joining $\sigma$ to $\sigma^\prime$.
For any $\sigma \in \cT$, we set
$$ \cT_\sigma = \big\{ \sigma^\prime \in \cT\, : \, \sigma \in \lgeo \rho , \sigma^\prime \rgeo     \big\} \; , $$
that is the subtree stemming from $\sigma$. We then set
$\Gamma (\cT_\sigma )= \sup_{\sigma^\prime \in \cT_\sigma} d(\sigma, \sigma^\prime)$ that is
the total height of $\cT_\sigma$. Next, for any $a, \varepsilon \in (0, \infty)$, we set
$$ \cT (a)= \big\{ \sigma \in \cT: d(\rho, \sigma )= a \big\} \quad {\rm and} \quad \cT (a, \varepsilon ) = \big\{ \sigma \in \cT (a):  \Gamma (\cT_\sigma ) > \varepsilon \big\} \; .$$
Since $\cT$ is a compact metric space $\cT (a, \varepsilon ) $ is a finite subset and we set
$Z_a (\varepsilon)= \# \cT (a, \varepsilon) $. Then,
$$ \cT(a, \varepsilon ) = \big\{ \sigma_1 , \ldots , \sigma_{Z_a (\varepsilon)} \big\} \; , $$
where, the $\sigma_i$ is the $i$-th point to be visited by $H$. 
For any $\eta \in (0, \infty)$, we set
$$  \cD_{a, \varepsilon, \eta} = \big\{ T_i \; ; \; 1 \leq i \leq Z_a (\varepsilon) \big\} \quad 
\textrm{where} \quad T_i =  \cT_{\sigma_i} \cap B(\sigma_i , \eta) .$$
The $T_i$s are the subtrees above level $a$ that are higher than $\varepsilon$ and 
that are truncated at height $\eta$. We simply call them the $(a, \varepsilon)$-{\it subtrees truncated at height $\eta$}.

  Recall from (\ref{heightform}) that $v(\varepsilon)=\baN (\Gamma (\cT) > \varepsilon)= ((\gamma-1)\varepsilon)^{-\frac{1}{\gamma-1}}$.
Recall from the subsection stating the branching property, that $\cG_a$ stands for the sigma-field generated by the height process below $a$. Recall  from (\ref{Nadef}) that $\baN_a= \baN (\, \cdot \, | \sup H >a )$. By the branching property at level $a$, 

\smallskip

\begin{itemize}

\item[($i$)] under $\baN_a$ and conditionally given $\cG_a$,
$Z_a (\varepsilon)$ is a Poisson random variable with parameter $v(\varepsilon)\langle \ell^a\rangle $. 

\end{itemize}

\smallskip

\noi 
Moreover, each $\cT_{\sigma_i}$ is coded by an excursion above level $a$ that is higher than $\varepsilon$. Therefore, 
\smallskip

\begin{itemize}

\item[($ii$)] conditionally given $Z_a (\varepsilon)$, the truncated subtrees $T_i$ are independent and distributed as $B(\rho, \eta)$ under $\baN_\varepsilon$.

\end{itemize}

\smallskip

\noi
Consequently, for any integer $k \geq 1$, for any measurable functions $F_1, \ldots , F_k : [0, \infty) \rightarrow [0, \infty)$, we have
\begin{equation}
\label{branchspecify}
 \baN_a \! \Big( \! \un_{\{ Z_a (\varepsilon)= k\}}\!\! \! \!\prod_{^{1\leq i\leq k}} \!\! \! F_i \big( \bm (T_i) \big)  \Big| \cG_a \Big) \!= \!
\frac{(v(\varepsilon)\langle \ell^a\rangle )^k}{k!}e^{-v(\varepsilon)\langle \ell^a\rangle }\!\! \! \prod_{^{1\leq i\leq k}} \!\! \! \baN_\varepsilon \big( F_i \big( \bm (B( \rho , \eta) \big) \big). \!\!
\end{equation}
The distribution of $\bm (B(\rho, \eta))$ under $\baN_\varepsilon $ plays an important role in the proofs and it is studied in Section \ref{tailsec}.

\medskip

The two following lemmas are used in the proofs of Theorem \ref{thin} and Theorem \ref{ThickB}.
We first show that any ball contains a reasonably large truncated subtree.
\begin{lemma}
\label{trianinball} Let  $r \in (0, 1/2)$. Let $n_r$ be the positive integer given by
$ 2^{-n_r} <r \leq 2^{-n_r+1}$. Let $\sigma \in \cT$ be such that $d(\rho, \sigma) \geq 2.2^{-n_r-2}$ and let $k \geq 1$ be the integer such that $(k+1)2^{-n_r-2} \leq d(\rho, \sigma) <
(k+2)2^{-n_r -2} $.

 Then, there exists a unique truncated subtree $T\in \cD_{k2^{-n_r-2}, 2^{-n_r-2}, 2^{-n_r-2}}$ such that $T \subset B(\sigma , r)$.
\end{lemma}
\noi
{\bf Proof:} to simplify notation, we set $a= k2^{-n_r-2}$ and  $\varepsilon = \eta =2^{-n_r-2}$. There
is a unique $i \in \{ 1, \ldots, Z_a (\varepsilon)\}$ such that $\sigma_i \in \lgeo \rho, \sigma\rgeo $ 
(and note that $d(\rho , \sigma_i )= a =  k2^{-n_r-2}$). Then,
$T_i = \cT_{\sigma_i} \cap B( \sigma_i , \eta) \subset B(\sigma, r)$. Indeed, for any 
$\sigma^\prime \in T_i$,  $ d( \sigma , \sigma^\prime) \leq d(\sigma, \sigma_i) + d(\sigma^\prime, \sigma_i)$. Since $\sigma^\prime \in B( \sigma_i , \eta)$, 
$d(\sigma^\prime, \sigma_i)\leq \eta =2^{-n_r-2}$, and since 
$\sigma_i \in \lgeo \rho, \sigma\rgeo $, we get  
$$ d(\sigma , \sigma_i) = d( \rho , \sigma) -d(\rho, \sigma_i) \leq (k+2)2^{-n_r -2} -a= 2. 2^{-n_r-2} . $$
Thus, $ d( \sigma , \sigma^\prime) \leq 3.2^{-n_r-2} < r$, which completes the proof. \cqfd

\medskip

Conversely, one proves that any ball is contained in a reasonably small truncated subtree.
\begin{lemma}
\label{ballintrian} Let  $r \in (0, 1/2)$. Let $n_r$ be the positive integer given by
$2^{-n_r} <r \leq 2^{-n_r+1}$. Let $\sigma \in \cT$ be such that $d(\rho, \sigma) \geq 2.2^{-n_r+1}$ and let $l \geq 1$, be the integer such that $(l+1)2^{-n_r+1} \leq d(\rho, \sigma) <
(l+2)2^{-n_r +1} $.

 Then, there exists a unique truncated subtree $T\in \cD_{l2^{-n_r+1}, 2^{-n_r+1}, 3.2^{-n_r+1}}$ such that $B(\sigma , r) \subset T$. 
\end{lemma}
\noi
{\bf Proof:} to simplify notation, we set $a= l2^{-n_r+1}$, $\varepsilon =2^{-n_r+1}$ and $\eta = 3 \varepsilon$. There is a unique
$i \in \{ 1, \ldots, Z_a (\varepsilon)\}$ such that $\sigma_i \in \lgeo \rho, \sigma\rgeo $ 
(and note that $d(\rho , \sigma_i )= a =  l2^{-n_r+1}$). We check that
$B(\sigma, r) \subset T_i = \cT_{\sigma_i} \cap B( \sigma_i , \eta) $. Indeed, let $\sigma^\prime \in B(\sigma, r)$ and let $\sigma \wedge \sigma^\prime$ be the branching point of $\sigma$ and $\sigma^\prime$. Namely, 
$ \lgeo \rho , \sigma \wedge \sigma^\prime \rgeo = \lgeo \rho , \sigma  \rgeo \cap \lgeo \rho , \sigma^\prime \rgeo $. We get 
$$ d( \rho , \sigma \wedge \sigma^\prime) = \min \big( d( \rho, \sigma) , d( \rho , \sigma^\prime ) \big) > d( \rho, \sigma ) -r \geq (l+1)2^{-n_r +1} -2^{-n_r +1} = l 2^{-n_r +1} \; .$$
Thus, $d( \rho , \sigma_i) \leq d (\rho, \sigma\wedge \sigma^\prime)$ and since $\sigma_i $ and $\sigma \wedge \sigma^\prime $ belongs to $ \lgeo \rho, \sigma\rgeo $, it implies that $\sigma^\prime \in \cT_{\sigma_i}$. Moreover, $ d(\sigma^\prime , \sigma_i) \leq d(\sigma^\prime , \sigma) + d(\sigma , \sigma_i) <r+ 2.2^{-n_r+1} \leq 3.2^{-n_r+1}= \eta$, 
which completes the proof. \cqfd

\section{Tail estimates.}
\label{tailsec}
\subsection{Tail of the distribution of $\bm (B(\rho, 1+c))$ under $\baN_1$.}
Recall that $v(1) = \baN (\sup H  >1)$ and that $\baN_1= \baN (\, \cdot \, | \sup H >1)$. Recall from 
(\ref{kappaexc}) the definition of $\kappa_a (\lambda, \mu)$ and Lemma \ref{kappatriangle} that asserts that for any $c\in [0, \infty)$, 
\begin{equation}
\label{recallPhi}
\baN_1 \big( e^{-\lambda \bm (B(\rho , 1+c)) }\big)= \Phi_{1, c} (\lambda) = v(1)^{-1} \big( \kappa_1  (\lambda, \infty ) -\kappa_{1+c} ( \lambda, 0) \big) \; .
\end{equation}
\begin{lemma}
\label{zerogammaun}
For any $\gamma \in (1, 2]$ and any $c \in [0, \infty)$, we get
$$ -\log \,  \baN_1 \! \big( \, \bm (B(\rho , 1+c) \, ) \leq y \, \big) \; \underset{^{y\rightarrow 0+}}{\sim} \; \Big(\frac{\gamma-1}{y} \Big)^{\gamma-1} \; .$$
\end{lemma}
\noi
{\bf Proof:} by  De Bruijn's Tauberian theorem, we only need to get an equivalent to $-\log \Phi_{1, c} (\lambda)$ when $\lambda$ goes to infinity. To that end, we use (\ref{recallPhi}) and we first get an estimate of $\kappa_1 (\lambda, \infty)$: we take $a= 1$ in (\ref{equakappa}) and we let $\mu$ go to infinity to obtain 
$$ \int_{\kappa_1(\lambda , \infty)}^\infty \frac{du}{u^\gamma -\lambda} = 1 \; .$$
We set $a(\lambda)= \log \big( \frac{\kappa_1 (\lambda, \infty)^\gamma}{ \lambda} -1\big)$ and we use the change of variable $y= \log (\lambda^{-1} u^\gamma -1)$ in the previous integral equation to get
$$ \gamma \lambda^{1-\frac{1}{\gamma}} = \int_{a(\lambda)}^\infty \frac{dy}{(1+e^y)^{1-\frac{1}{\gamma}}} \; .$$
Note that $a(\lambda)$ decreases to $ -\infty$ when $\lambda$ goes to $\infty$. Thus, there exists $\lambda_0 \in (0, \infty)$ such that $a(\lambda) <0$, for any $\lambda >\lambda_0$.
Next, observe that 
\begin{equation}
\label{Qzero}
Q_0(\lambda) \!\! : = \int_0^\infty \!\!  \frac{dy}{(1+e^y)^{1-\frac{1}{\gamma}}} - \int_{a(\lambda)}^0
\!\! \Big(1-\frac{1}{(1+ e^y)^{1-\frac{1}{\gamma}  }} \Big) dy \; \underset{^{\lambda \rightarrow \infty}}{-\!\!\! -\!\!\! \longrightarrow}  Q_0(\infty) \in \bR \, ,
\end{equation}
and that $\gamma \lambda^{1-\frac{1}{\gamma}} = Q_0(\lambda) -a(\lambda) $, for any $\lambda >\lambda_0$. Namely,
\begin{equation}
\label{kappainfty}
\kappa_1 (\lambda , \infty) = \lambda^{\frac{1}{\gamma}} \Big( 1+
\exp \big(Q_0 (\lambda) -\gamma \lambda^{\frac{\gamma-1}{\gamma}} \big)  \Big)^{\frac{1}{\gamma}} \; .
\end{equation}

Similarly, we get an estimate for $\kappa_{1+c} (\lambda, 0) $: we take $a= 1+c$ and $\mu= 0$ in (\ref{equakappa}):  
 $$ \int_0^{\kappa_{1+c}(\lambda , 0)} \frac{du}{\lambda - u^\gamma } = 1+c  \; .$$
We set $b(\lambda)= -\log \big( 1-\frac{\kappa_{1+c} (\lambda, 0)^\gamma}{ \lambda}\big)$
and we  take $y= -\log (1-\lambda^{-1} u^\gamma )$ to get
$$ (1+c)\gamma \lambda^{1-\frac{1}{\gamma}} = \int_{0}^{b(\lambda)} \frac{dy}{(1-e^{-y})^{1-\frac{1}{\gamma}}}=  b(\lambda) + Q_1(\lambda) \; , $$
where
\begin{equation}
\label{Qun}
Q_1(\lambda): = \int_0^{b(\lambda)}  \!\! \Big( \;  \frac{1}{\; (1- e^{-y})^{1-\frac{1}{\gamma}  }} -1 \Big) dy \; \underset{^{\lambda \rightarrow \infty}}{-\!\!\! -\!\!\! \longrightarrow}  Q_1(\infty) \in [0, \infty ) \, .
\end{equation}
Thus,
\begin{equation}
\label{kappaainfty}
\kappa_{1+c} (\lambda , 0) = \lambda^{\frac{1}{\gamma}} \Big( 1-
\exp \big(Q_1 (\lambda) - (1+c) \gamma \lambda^{\frac{\gamma-1}{\gamma}} \big)  \Big)^{\frac{1}{\gamma}} \; .
\end{equation}
By (\ref{recallPhi}), (\ref{Qzero}), (\ref{kappainfty}), (\ref{Qun}) and (\ref{kappaainfty}), we get
$$ - \log \baN_1 \big( e^{-\lambda \bm (B(\rho , 1+c)) }\big) \; \underset{^{\lambda \rightarrow \infty}}{\sim}
\;   \gamma \lambda^{\frac{\gamma-1}{\gamma}} \; , $$
and De Bruijn's Tauberian theorem entails the desired result  (see Theorem 4.12.9 page 254 in \cite{BiGoTe}). \cqfd

\medskip

 In the Brownian case, computations are explicit: we easily derive from (\ref{equakappa}) that for any $a, \lambda , \mu\in [0, \infty)$ such that $\sqrt{\lambda} \neq \mu$,
\begin{equation}
 \label{explicitkappa}
\kappa_a (\lambda , \mu)= \sqrt{\lambda} \cdot  \frac{e^{a\sqrt{\lambda}} (\sqrt{\lambda} +\mu) -
e^{-a\sqrt{\lambda}}(\sqrt{\lambda} -\mu)  }{e^{a\sqrt{\lambda}} (\sqrt{\lambda} +\mu) +
e^{-a\sqrt{\lambda}}(\sqrt{\lambda} -\mu)  } \; .
\end{equation}
Recall that $\coth (x)= (e^x+e^{-x})/(e^x-e^{-x})= 1/ \tanh (x)$, and note that $v(1)= 1$. 
Thus, (\ref{recallPhi}) implies that 
$$  \baN_1 \big( e^{-\lambda \bm (B(\rho , 1+c)) }\big)= \sqrt{\lambda} \big( \coth (\sqrt{\lambda})- \tanh ((1+c)\sqrt{\lambda}) \big) \; .$$
We next use the well-known formulas
\begin{equation}
\label{thanhdevel}
 x\tanh (x)=  \sum_{n \geq 0} \frac{2x^2}{x^2+ \frac{\pi^2}{4}(2n+1)^2} \quad {\rm and } \quad x\coth (x)= 1+ \sum_{n \geq 1} \frac{2x^2}{x^2+ \pi^2 n^2} \; .
\end{equation}
Therefore,
\begin{eqnarray*}
\int_0^\infty \!\! \baN_1 \big( \bm (B(\rho , 1+c)) \geq  y \big) e^{-\lambda y}\, dy  & =& \lambda^{-1} \big( 1-  \baN_1 \big( e^{-\lambda \bm (B(\rho , 1+c)) }\big) \, \big) \\
&= & \frac{_2}{^{1+c}} \sum_{n \geq 0} \frac{1}{\lambda + \big(\frac{\pi (2n+1)}{2(1+c)}\big)^2} -2 \sum_{n \geq 1} \frac{1}{\lambda + \pi^2 n^2}.
\end{eqnarray*}
This easily implies the following.
\begin{lemma}
\label{inverBrown} Consider the Brownian case: $\gamma=2$. Then, for any $y \in [0, \infty)$,
\begin{eqnarray*}
\baN_1 \big( \bm (B(\rho , 1+c)) \geq  y \big) &= &  \frac{_{2}}{^{1+c}} \sum_{n \geq 0}
\exp \big(-\frac{_{\pi^2(2n+1)^2}}{^{4(1+c)^2}} \, y  \big) - 2 \sum_{n \geq 1} \exp (-\pi^2 n^2 \, y) \\
& \underset{^{y \rightarrow \infty}}{\sim} &  \frac{_2}{^{1+c}} \,  \exp \big(-\frac{_{\pi^2}}{^{4(1+c)^2}} \, y  \big) \; .
\end{eqnarray*}
\end{lemma}
This result shall be used in the proof of Therem \ref{ThickB}.

\subsection{Tail of the distribution of $M_*$.}
Recall from (\ref{Mstarun}) that $M_*= M^*_1 (1)$. We set
$ \cL (\lambda) = \bE \big[ \exp (-\lambda M_*) \big]$ and from (\ref{QLapcour}), we get  
\begin{equation}
\label{MstarLapl}
 \cL (\lambda) = \bE \big[ \exp (-\lambda M_*) \big]= 1- \frac{\kappa_1(\lambda, 0)^\gamma}{\lambda} \, , \quad \lambda \in [0, \infty) \; .
\end{equation}
The following lemma provides an equivalent of the tail at $0+$ of the distribution of $M_*$ that is used in
the proof of Proposition \ref{Constant}.
\begin{lemma}
\label{masszero} For any $\gamma \in (1, 2]$ we have the following estimate.
$$ \lim_{y \rightarrow 0 +} \; y^{-\frac{\gamma -1}{2}}
\exp \!\Big( \, \frac{1}{\; y^{\gamma -1}} \, \Big) \;
 \bP \big(  M_* \leq  (\gamma \!- \!1) \,  y \big)  \; = \; e^{C_\gamma } \sqrt{\frac{\gamma (\gamma-1)}{2 \pi}}   \; .$$
where $C_\gamma$ is a constant given by
\begin{equation}
\label{Cdeff}
 C_\gamma = \int_0^1  u^{-1} \big(   (1-u)^{-\frac{\gamma-1}{\gamma}} -1 \big) du
  = \sum_{n \geq 1}
\frac{1}{n} \left| \binom{\frac{1-\gamma}{\gamma}}{n}\right|  \; ,
\end{equation}
\end{lemma}
\noi
{\bf Proof:} first observe that for any $u \in [0, 1)$, we have
\begin{equation}
\label{keypower}
(1-u)^{-\frac{\gamma-1}{\gamma} } =  \sum_{n \geq 0}  (-1)^n \binom{\frac{1-\gamma}{\gamma}}{n} u^n =   1 + \sum_{n \geq 1} \left|
\binom{\frac{1-\gamma}{\gamma}}{n} \right|  u^n \; .
\end{equation}
This easily entails the second equality in (\ref{Cdeff}). For any $y\in [0, 1)$, we set
$$ F(y) = \int_{y}^1 \frac{du}{u(1-u)^{\frac{\gamma-1}{\gamma} } } \; $$
By (\ref{equakappa}) and by a simple change of variable, $F( \cL (\lambda))= \gamma \lambda^{\frac{\gamma-1}{\gamma}} $. Note that
\begin{eqnarray*}
F(y)& =& \int_y^1 u^{-1} du   \; + \int_0^1 \!\!\!  u^{-1} \big(   (1-u)^{-\frac{\gamma-1}{\gamma}} -1 \big) du  \;  -\int_0^y  \!\!\! u^{-1} \big(   (1-u)^{-\frac{\gamma-1}{\gamma}} -1 \big) du  \\
& =& -\log y  + C_\gamma - h(y) \; .
\end{eqnarray*}
Here $h: [0, 1] \rightarrow [0, \infty) $ is increasing, $h(0)=0$,
$h(1)= C_\gamma$, and $ h(y)= \sum_{n \geq 1} a_n y^n $, where for any $n \geq 1$,
\begin{equation}
\label{coeffh}
 \frac{h^{(n)} (0)}{n!} = \frac{1}{n}\left|
\binom{\frac{1-\gamma}{\gamma}}{n} \right|= \frac{1}{n} \prod_{k=1}^n \left( 1- \frac{1}{\gamma  \, k}\right) \in (0, 1)   \; .
\end{equation}
Thus,
\begin{equation}
\label{integcatr}
\cL (\lambda) = e^{C_\gamma} \exp \! \big( \!- \! \gamma \lambda^{\frac{\gamma-1}{\gamma}} \! \big) \, \exp \left( -h(\cL (\lambda) )\right) \; .
\end{equation}
We next use Fubini for sums of nonnegative real numbers to prove that for any $y\in [0, 1]$,
and any integer $m \geq 1$,
\begin{equation}
\label{hmpower}
h(y)^m \! = \!\!  \sum_{n \geq m} y^n \;  \cdot  \! \! \! \!  \!\!  \!\!\!\! \!\! \sum_{\; \; \; \substack{ q_1, \ldots , q_m \geq 1\\ q_1+ \ldots + q_m=n  }}    \!\!\!\!\! \!  a_{q_1} \ldots a_{q_m} =
\sum_{n \geq m} y^n \;  \cdot   \! \! \! \!  \! \! \! \!  \! \! \! \!  \! \! \! \!  \sum_{\; \; \; \; \substack{ p_1+ \ldots + p_n= m \\ p_1+ 2p_2+ \ldots + np_n=n  }}  \!\!\!\!      \frac{m!}{p_1! \ldots p_n!} a_{1}^{p_1} \ldots a_{n}^{p_n} .
\end{equation}
Thus, for any $y \in [0, 1]$, $ \exp (h(y))=  1 + \sum_{n \geq 1} d_n y^n$,
$ d_n =\sum a_1^{p_1} \ldots  a_n^{p_n}/ (p_1! \ldots p_n!)$, where the sum is over all the $p_1,  \ldots , p_n \geq 0$ such that $p_1+ 2p_2 + \ldots + np_n= n$.
Standard arguments on analytic functions  imply that there exists $r_1 >0$ such that
$$ \exp (-h(y))=  1+ \sum_{n \geq 1} c_n y^n  \; , \quad y\in [0, r_1)  $$
and (\ref{hmpower}) easily entails that for any $n\geq 1$,
$$ c_n= \sum_{\substack{ p_1,  \ldots , p_n \geq 0 \\ p_1+ 2p_2 + \ldots + np_n= n}} \frac{(-1)^{p_1+ \ldots + p_n}}{p_1! \ldots p_n! }a_1^{p_1} \ldots  a_n^{p_n} \; .$$
Consequently, $|c_n| \leq d_n $, and $r_1 \geq 1$, which implies that
\begin{equation}
\label{boundsecu}
\sum_{n \geq 0} | c_n| \leq \exp (h(1))= e^{C_\gamma} \; .
\end{equation}
The previous arguments and (\ref{integcatr}) imply
\begin{equation}
\label{integfive}
\cL (\lambda) = e^{C_\gamma} \exp \big( \!- \! \gamma \lambda^{\frac{\gamma-1}{\gamma}} \big) +
\sum_{n \geq 1} e^{C_\gamma}c_n  \exp \big( \!- \! \gamma \lambda^{\frac{\gamma-1}{\gamma}} \big) \cL (\lambda)^n \; , \quad \lambda \in [0, \infty) \; .
\end{equation}
We next introduce a non-negative random variable $S $ defined on $(\Omega , \cF , \bP)$,  that has a stable distribution whose Laplace transform is given by
$$ \bE \big[ \exp \big( \! -\! \lambda S \big)  \big] = \exp \big( \!- \! \gamma \lambda^{\frac{\gamma -1}{\gamma}} 
\big) \; .$$
We use the following standard tail estimate due
to Skorohod \cite{Sko54} (see also Example 4.1 in \cite{JaPr87}).
\begin{equation}
\label{Zolo}
 \lim_{y \rightarrow 0 +} \; y^{-\frac{\gamma -1}{2}}
\exp \!\Big( \, \frac{1}{\; y^{\gamma -1}} \, \Big) \;
\bP \big(  S  \leq  (\gamma \!- \!1) \,  y \big)
\; = \; \sqrt{\frac{\gamma (\gamma-1)}{2 \pi}} \; .
\end{equation}
We denote by $q$ the density of $S$ and by $\mu$ the distribution of $M_*$. The bound (\ref{boundsecu}) implies that
$$ R (dx) := e^{C_\gamma} q(x) dx + \sum_{n \geq 1} e^{C_\gamma}c_n (q\ast \mu^{\ast n } ) (dx) $$
is a Borel signed measure on $[0, \infty)$ whose total variation is bounded by $2 \exp (2C_\gamma)$.
Moreover (\ref{integfive}), implies that $\int_{[0, \infty)} e^{-\lambda x} R(dx) = \cL(\lambda)$, $\lambda \in [0, \infty)$.
Standard arguments on Laplace transform imply that $R= \mu $.
Denote by $(Y_n )_{n \geq 1}$, a sequence of i.i.d.$\;$copies of $M_*$ that are also independent of $S$. Since $R= \mu$, for any $y \in [0, \infty)$, 
$$ \bP( M_* \leq (\gamma-1) y) = e^{C_\gamma} \bP (S  \leq (\gamma-1) y ) + \sum_{n \geq 1} e^{C_\gamma}
c_n \bP (S + Y_1 + \ldots + Y_n \leq (\gamma-1) y)  .$$
The obvious bound
$$\bP \big(  S + Y_1 + \ldots + Y_n \leq (\gamma-1) y \,  \big) \, \leq  \, \bP \big( \, S \leq (\gamma-1) y \big) \, \bP \big( \, M_* \leq (\gamma-1) y \big)^n \;  $$
entails
$$\left| \frac{ \bP (M_* \leq (\gamma-1) y ) }{ e^{C_\gamma} \bP (S \leq (\gamma-1) y ) } -1 \right|  \quad \leq \quad  
\sum_{n\geq 1} |c_n| \; \bP \big(  M_* \leq (\gamma-1) y \big)^n  \quad \underset{^{y\rightarrow 0+}}{-\!\!\!-\!\!\!-\!\!\!-\!\!\! \longrightarrow} \; 0 \; , $$
which entails the desired result thanks to (\ref{Zolo}). \cqfd 

\bigskip

In the Brownian case, the computations are explicit. By (\ref{explicitkappa}), we get
$$ \int_0^\infty e^{-\lambda y} \bP (M_* \geq y) = \lambda^{-1} \big( 1-\cL (\lambda) \big) =
\lambda^{-1} \tanh^2 (\sqrt{\lambda} ) \; .$$
Observe that $\tanh^2 (x)= \tanh^\prime (0)-\tanh^\prime (x)$. If we set 
$a_n = \pi^2 (2n+1)^2/ 4$, then (\ref{thanhdevel}) implies 
$$ \tanh^2 (x)= \sum_{n \geq 0} \; \frac{2x^2}{a_n (x^2+a_n)}  +\frac{4x^2}{(x^2+a_n)^2} .$$
Thus
$$ \int_0^\infty e^{-\lambda y} \bP (M_* \geq y) = \sum_{n \geq 0} \frac{2}{a_n (\lambda +a_n)} +
\frac{4}{(\lambda +a_n)^2} \; , $$
which easily implies the following lemma.
\begin{lemma}
\label{inverBrownMstar} Consider the Brownian case: $\gamma=2$. Then, for any $y \in [0, \infty)$,
\begin{eqnarray*}
\bP \big( M_* \geq  y \big) &= &  \sum_{n \geq 0}  4 \big( \; \frac{_2}{^{\pi^2 (2n+1)^2} }  + y    \big) \exp \big( \! -\! \frac{_{\pi^2}}{^4}  (2n+1)^2 y \big)  \\
& \underset{^{y \rightarrow \infty}}{\sim} &  4y   \exp \!\big(\!-\!\frac{_{\pi^2}}{^{4}} \, y  \big) \; .
\end{eqnarray*}
\end{lemma}
This lemma shall be used in the proof of Proposition \ref{constant}.

\section{Proofs.}
\label{proofsec}
\subsection{Proof of Theorem \ref{thin}.}
\label{thinsec}
We fix $\gamma \in (1, 2]$ and we consider the $\gamma$-stable tree $(\cT, d)$ with root $\rho$ coded by the $\gamma$-height process $H$ under its excursion measure $\baN$, as defined in the introduction. Recall from Section \ref{Balltrianglesec} the definition of the $(a, \varepsilon )$-subtrees truncated at height $\eta$, whose set is denoted by
$\cD_{a, \varepsilon, \eta}= \{ T_i \,  ; \,  1\leq i\leq Z_a (\varepsilon) \}$.
Recall that
$$f_\gamma (r)= \frac{r^{\frac{\gamma}{\gamma-1} }}{ (\log 1/r)^{\frac{1}{\gamma-1}}} \, , \quad
r\in (0,  1) \; .$$

\medskip

\noi
{\bf Lower Bound.} Let fix a positive integer $R_0$ and a real number $\alpha \in (0, \infty)$, that is specified further. For any integer $n\geq 4$, we set
$$ V(n)= \un_{\{ \bm (B(\rho , 2^{-n}) \leq   \alpha f_\gamma (2^{-n})\}  }  \; +\!\! \! \sum_{1\leq k < R_02^n}
\!\! \!\! \# \big\{ T_i \in \cD_{k2^{-n}, 2^{-n}, 2^{-n}} \; : \;
\bm (T_i)  \leq \alpha f_\gamma (2^{-n})  \big\}  .$$
We first prove that
\begin{equation}
\label{Vbascule}
 V(n)= 0 \quad  \Longrightarrow \quad  \inf_{\sigma \in B(\rho , R_0) } \bm \big( B(\sigma , 2^{-n+3}) \big) \;  > \, \alpha \,  f_\gamma (2^{-n}) \; .
\end{equation}
Indeed, we apply Lemma \ref{trianinball} with $r= 2^{-n+3}$. Thus, $n_r= n-2$ and we have
$2^{-n_r} < r = 2^{-n_r +1}$.  Let $\sigma \in \cT$ be such that $d(\rho , \sigma) \leq R_0$. We first consider the case where $ d(\rho , \sigma ) \geq
2.2^{-n_r-2} = 2^{-n +1}=r/4$. Let $k \in \bN$ be such that $ (k+1)2^{-n_r-2} \leq d(\rho, \sigma ) < (k+2)2^{-n_r -2}$. Observe that $1 \leq k < R_0 2^n$.  Lemma
\ref{trianinball} implies that there exists a truncated subtree
$T\in \cD_{k2^{-n}, 2^{-n} , 2^{-n}}$ such that
$T \subset B( \sigma , 2^{-n+3} )$. Consequently, if $V(n)= 0$, then $\alpha f_\gamma (2^{-n} ) <
\bm (B( \sigma , 2^{-n+3} ))$.

We next consider the case where $d(\rho , \sigma) <  2^{-n+1}$. Then it is easy to see that
$B( \rho , 2^{-n})\subset  B( \sigma , 3. 2^{-n}) \subset B( \sigma , 2^{-n+3} )$. Thus, if $V(n)= 0$,
$\bm (B(\sigma , 2^{-n +3})) \geq \bm (B(\rho , 2^{-n})  >  \alpha f_\gamma (2^{-n})$,
which completes the proof of (\ref{Vbascule}).

\medskip

  We next claim that it is possible to find $\alpha$ such that
\begin{equation}
\label{VsumBC}
\sum_{n \geq 4} \baN \big( V(n) \un_{\{ \Gamma (\cT ) > 2^{-n} \} } \big) < \infty \; .
 \end{equation}
We first set $x_n= \baN \big( V(n) \un_{\{ \Gamma (\cT ) > 2^{-n} \} } \big) $ and
\begin{eqnarray*}
y_n & = & \baN \big(  \; \bm (B(\rho , 2^{-n} )) \! \leq  \!  \alpha f_{\gamma} (2^{-n})  \quad  {\rm and } \quad  \Gamma (\cT) \! > \! 
2^{-n} \;  \big) \\
 &= & v(2^{-n} )\baN_{2^{-n}} \big( \; \bm (B(\rho , 2^{-n} )) \! \leq \! \alpha f_{\gamma} (2^{-n}) \;   \big) \;.
\end{eqnarray*}
To simplify notation, we also set $Z(k,n)=  Z_{k2^{-n} } (2^{-n})$. We  get the following.
\begin{eqnarray*}
x_n & \leq & y_n  +
 \sum_{1\leq k<R_02^n}  \baN \Big( \sum_{^{1\leq i \leq Z(k,n)  } } \!\!\!\! \un_{\{ \bm(T_i) \leq \alpha f_{\gamma} (2^{-n})  \} }  \Big) \\
& \leq & y_n   +
 \sum_{1\leq k<R_02^n } v(k2^{-n}) \baN_{k2^{-n}}  \Big( \sum_{^{1\leq i \leq Z(k,n)  } } \!\!\!\!
 \un_{\{ \bm (T_i)\leq \alpha f_{\gamma} (2^{-n})  \} }  \Big)
\end{eqnarray*}
Recall from the definition of the branching property that $\cG_{k2^{-n}}$ stands for 
the sigma-field generated by the tree below $k2^{-n}$. Recall from Section \ref{Balltrianglesec} that conditionnally given $\cG_{k2^{-n}}$, $Z(k,n)$ is a Poisson random variable with parameter $v(2^{-n}) \langle \ell^{k2^{-n}} \rangle $ and by 
 (\ref{branchspecify}), we get 
\begin{eqnarray*}
 \baN_{k2^{-n}} \! \Big( \sum_{^{1\leq i \leq Z(k,n)  } } \!\!\!\!
 \un_{\{ \bm(T_i) \leq \alpha f_{\gamma} (2^{-n})  \} }  \Big)  \!\!\!\!  \!\!  & = &  \!\!\!\!    \baN_{k2^{-n}} \big(Z(k, n) \big)  \,
\baN_{2^{-n}} \! \big( \bm (B(\rho , 2^{-n} )) \!  \leq \!  \alpha f_{\gamma} (2^{-n})  \big)  \\
& =&   \!\!\!\!    \baN_{k2^{-n}} \big( \langle \ell^{k2^{-n}} \rangle  \big) \, y_n \; .
\end{eqnarray*}
Consequently, $ x_n \;    \leq \;  y_n  +  \sum_{1\leq  k<R_02^n} v(k2^{-n}) \baN_{k2^{-n}} \big( \langle \ell^{k2^{-n}} \rangle  \big)  \, y_n$. 
We next use (\ref{Nabranchmean}) that implies $ v(k2^{-n}) \baN_{k2^{-n}} \big( \langle \ell^{k2^{-n}} \rangle  \big) = 1$. Thus, we get $x_n \leq  R_0 2^{n} y_n $. Namely,
\begin{equation}
\label{namelyun}
\baN \big( V(n) \un_{\{ \Gamma (\cT ) > 2^{-n} \} } \big)  \leq R_0\,  2^n v(2^{-n} )\; \baN_{2^{-n}} \big( \bm (B(\rho , 2^{-n} )) \leq  \alpha f_{\gamma} (2^{-n})  \big)  \; .
\end{equation}
By (\ref{heightform}), $2^n v(2^{-n})= (\gamma-1)^{-\frac{1}{\gamma-1}} 2^{\frac{ \gamma n  }{\gamma-1}}$ and the scaling property (\ref{trianscale}) implies
$$ \baN_{2^{-n}} \big( \bm (B(\rho , 2^{-n} )) \leq  \alpha f_{\gamma} (2^{-n})  \big)  =
\baN_{1} \big( \bm (B(\rho , 1 )) \leq  \alpha (\log 2^{n})^{-\frac{1}{\gamma-1}}  \big) \; .$$
By Lemma \ref{zerogammaun}, there is a constant $q  \in (0, \infty)$
that only depends on $\gamma$, such that for any $n \geq 4 $,
$$   \baN_{1} \big( \bm (B(\rho , 1 )) \leq  \alpha (\log 2^{n})^{-\frac{1}{\gamma-1}}  \big) \leq
\exp \big( - \alpha^{-(\gamma-1)} q \log 2^n \big)   \; .$$
This inequality combined with (\ref{namelyun}) entails that  
for any $n \geq 4$, 
$$ \baN \big( V(n) \un_{\{ \Gamma (\cT ) > 2^{-n} \} } \big) \leq R_0(\gamma-1)^{-\frac{1}{\gamma-1}}  \exp \big(- \big(  \alpha^{-(\gamma-1)} q  \log 2- \frac{_{\gamma}}{^{\gamma -1}}\log 2 \big) \,  n   \big) \, ,
 $$
which implies (\ref{VsumBC}) if $\alpha < \big(\frac{(\gamma-1)q  }{ \gamma }
\big)^{\frac{1}{\gamma -1}} $.

\medskip

Since $V(n) \in \bN$, (\ref{VsumBC}) implies that $\baN$-a.e.$\;$for any sufficiently large $n$, we have $V(n)\un_{\{ \Gamma (\cT) >2^{-n}\}}= 0$. Since $\baN (\Gamma (\cT) =0)=0$, it implies that $\baN$-a.e.$\;$for any sufficiently large $n$,
$V(n)= 0$. We next use (\ref{Vbascule}), to get the following: there exists $\alpha_\gamma \in (0, \infty)$,  that only depends on $\gamma$, such that for any positive integer $R_0$,
$$ \textrm{$\baN$-a.e.$\; \; \; \exists n_0\in \bN$ s.t. $ \forall n \geq n_0$,}  \quad  \inf_{\sigma \in B(\rho , R_0) } \bm \big( B(\sigma , 2^{-n+3}) \big)  \geq \alpha_\gamma \, f_\gamma (2^{-n} ) \; .$$
Note that $\baN $-a.e.$\;$there exists $R_0$ such that $B( \rho,  R_0)= \cT$. Since $\alpha_\gamma$ does not depend on $R_0$, it entails 
\begin{equation}
\label{lowdyad}
\textrm{$\baN$-a.e.}  \qquad \liminf_{n \rightarrow \infty} \; \frac{1}{f_\gamma (2^{-n})}  \;
\inf_{\sigma \in \cT} \bm \big( B(\sigma , 2^{-n+3}) \big)  \geq \alpha_\gamma  \; .
\end{equation}

\medskip

\noi
{\bf Upper Bound.} Let $R_0$ be a positive integer and let $h_0 \in (0, \infty)$. We also fix $\beta \in (0, \infty)$, that is specified further. We introduce the following event
$$ A_n = \big\{ \Gamma (\cT) > h_0  \big\} \cap \big\{  \inf_{^{\sigma \in B(\rho , R_0)}}
\bm (B (\rho , 2^{-n}) )   > \beta  f_\gamma (2^{-n} ) \,  \big\} \; . $$
We assume that $n \geq 4$, and that $h_0 > 2^{-n+1}$. Let $l \geq 1$ be such that
$(l+1)2^{-n} \leq h_0 < (l+2)2^{-n}$. We argue on the event $A_n$:  let $\sigma \in \cT$ be such that $d(\rho, \sigma )= h_0$; we apply Lemma \ref{ballintrian} with $r= 2^{-n}$; thus
$n_r= n+1$ (namely, $2^{-n-1}= 2^{-n_r} < r= 2^{-n_r+1}= 2^{-n}$); it implies that there exists a truncated subtree 
$T \in \cD_{l2^{-n} , 2^{-n} , 3.2^{-n}} $ such that $B(\sigma, 2^{-n}) \subset T$. This proves 
\begin{equation}
\label{AsubA}
A_n \subset B_n := \big\{ \Gamma (\cT) > l2^{-n}  \big\} \cap \big\{  \forall T_i \in  \cD_{l2^{-n} , 2^{-n} , 3.2^{-n}} : \bm (T_i) > \beta f_\gamma (2^{-n} )  \big\} \; .
\end{equation}
To simplify notation we set $Z(l,n)= Z_{l2^{-n}} (2^{-n} ) $. We then use (\ref{branchspecify}) to get
\begin{eqnarray*}
\baN (B_n) & = & \baN (\Gamma (\cT) >l2^{-n} ) \, \baN_{l2^{-n}} \Big( \prod_{^{1\leq i\leq Z(l,n)}}
\un_{\{ \bm (T_i) > \beta f_\gamma (2^{-n} ) \}}  \Big) \\
& =& v(l2^{-n}) \,  \baN_{l2^{-n}} \Big( \, \baN_{2^{-n}} \! \big( \,  \bm (B(\rho, 3.2^{-n} ) ) \!> \!  \beta
f_\gamma (2^{-n}) \,  \big)^{Z(l,n)} \Big) \; .
\end{eqnarray*}
Recall that under $\baN_{l2^{-n}}$, conditionally given $\cG_{l2^{-n}}$, the random variable 
$Z(l,n)$ has a Poisson distribution with mean $v(2^{-n}) \langle \ell^{l2^{-n}} \rangle $. We then get
\begin{equation}
\label{EnnBestep}
\baN (B_n)= v(l2^{-n}) \;  \baN_{l2^{-n}} \left( \exp \left( -z_n \langle \ell^{l2^{-n}} \rangle  \right) \right) \; ,
\end{equation}
where we have set
$$z_n= v(2^{-n})\;  \baN_{2^{-n}}  \big( \,  \bm (B(\rho, 3.2^{-n} ) ) \leq  \beta
f_\gamma (2^{-n}) \,  \big)\; .$$
By (\ref{heightform}) and the scaling property (\ref{kappascale}), we get 
$$ z_n= (\gamma-1)^{-\frac{1}{\gamma-1}} 2^{\frac{n}{\gamma-1}} \baN_1 \big( \bm(B(\rho, 3) \leq \beta (\log 2^n)^{-\frac{1}{\gamma- 1}} \,   \big) \; .$$
By Lemma \ref{zerogammaun}, there exists $q, q^\prime \in (0, \infty)$ that only depend on $\gamma$ such that for any $n \geq 4$,
$$ z_n \;  \geq \; w_n := q^\prime \exp \left(   \left(  \frac{_1}{^{\gamma-1}} \log 2 - \beta^{-(\gamma -1) } 
q  \log 2 \right)  n \right)  \; .$$
We fix $\beta > \big( (\gamma-1)q \big)^{\frac{1}{\gamma-1}}$ so that $\theta_\gamma :=  \frac{1}{{\gamma-1}} \log 2 - \beta^{-(\gamma -1) }q  \log 2 >0$. Thus, $w_n= q^\prime 
 \exp (\theta_\gamma n) \rightarrow \infty$, when $n \rightarrow \infty$.
By  (\ref{EnnBestep}), we get
$$ \baN (B_n) \leq  v(l2^{-n}) \,  \baN_{l2^{-n}} \big( \exp \big(-w_n \langle \ell^{l2^{-n}} \rangle  \big) \big) \; . $$
Recall that $h_0 > 2^{-n+1}$, which implies that $l \geq 1$. Thus $l2^{-n} \geq h_0/3$. Since $v$ decreases, $ v(l2^{-n}) \leq v(h_0/3)$. Next, recall (\ref{Nabranch}) and observe that $a\mapsto \baN_a (\exp (-\mu \langle  \ell^a\rangle  ))$ is decreasing. Thus,
\begin{equation}
\label{EnnBestepdeuz}
 \baN (B_n) \leq  v(h_0/3) \,  \baN_{h_0/3} \big( \exp \big(-w_n \langle \ell^{h_0/3} \rangle  \big) \big) \; .
\end{equation}
Since $\lim_{n \rightarrow \infty} w_n = \infty$, we easily derive from (\ref{Nabranch}) with $a= h_0/3$ that
$$ \baN_{h_0/3} \big( \exp \big(-w_n \langle \ell^{h_0/3} \rangle  \big) \big)  \ \underset{^{n \rightarrow \infty}}{\sim} \; \frac{_{3}}{^{(\gamma-1)^2(q^\prime)^{\gamma-1} h_0}} \, \exp \!\big(\! -\! (\gamma\!-\!1) \theta_\gamma  n \big) \; . $$

  Thus, and (\ref{AsubA}) and (\ref{EnnBestepdeuz}) immediately entail $ \sum_{n \geq 4} \baN (A_n) < \infty$. By Borel-Cantelli, there exists $\beta_\gamma \in (0, \infty)$, that only depends on $\gamma$, such that for any $h_0, R_0 \in (0, \infty)$, $\baN$-a.e.$\;$for any sufficiently large $n$, $\un_{A_n} = 0$. Since $\baN$-a.e.$ \, 0<\Gamma (\cT) <\infty$, one gets
\begin{equation}
\label{uppdyad}
\textrm{$\baN$-a.e.}  \qquad \limsup_{n \rightarrow \infty} \;  \frac{1}{f_\gamma (2^{-n})}  \;
\inf_{\sigma \in \cT} \bm \big( B(\sigma , 2^{-n}) \big)  \leq  \beta_\gamma  \; .
\end{equation}
Since $f_\gamma$ is regularly varying at $0$, (\ref{lowdyad}) and (\ref{uppdyad}) entail Theorem \ref{thin}. \cqfd

\subsection{Proof of Theorem \ref{ThickB}.}
\label{ThickBsec}
The proof is close to that of Theorem \ref{thin}. Here, we fix $\gamma= 2$ and we recall that 
$ f(r)= r^2 \log 1/r $, $r\in (0, 1)$.

\medskip

\noi
{\bf Upper Bound.} We fix a positive integer $R_0$ and a real number $\alpha \in (0, \infty)$ that is specified further. For any integer $n\geq 4$, we set
$$ W(n)= \un_{\{ \bm (B(\rho , 3.2^{-n}) \geq   \alpha f (2^{-n})\}  }  \; +\!\! \! \sum_{1\leq l < R_02^n}
\!\! \!\! \# \big\{ T_i \in \cD_{k2^{-n}, 2^{-n}, 3.2^{-n}} \; : \;
\bm (T_i)  \geq \alpha f (2^{-n})  \big\}  .$$
Arguing as previously, we apply Lemma \ref{ballintrian} with $r=2^{-n}$ to prove that
\begin{equation}
\label{Wbascule}
 W(n)= 0 \quad  \Longrightarrow \quad  \sup_{\sigma \in B(\rho , R_0) } \bm \big( B(\sigma , 2^{-n}) \big) \;  \leq  \, \alpha \,  f (2^{-n}) \; .
\end{equation}

 We next claim that it is possible to find $\alpha$ such that
\begin{equation}
\label{WsumBC}
\sum_{n \geq 4} \baN \big( W(n) \un_{\{ \Gamma (\cT ) > 2^{-n} \} } \big) < \infty \; .
 \end{equation}
To simplify notation we first set $x_n= \baN \big( W(n) \un_{\{ \Gamma (\cT ) > 2^{-n} \} } \big) $ and
\begin{eqnarray*}
y_n & = & \baN \big( \, \bm (B(\rho , 3. 2^{-n} )) \! \geq \!  \alpha f (2^{-n})  \quad  {\rm and } \quad  \Gamma (\cT) \! > \! 
2^{-n} \, \big) \\
 &= & v(2^{-n} )\baN_{2^{-n}} \big( \bm (B(\rho , 3. 2^{-n} )) \geq \alpha f (2^{-n})  \big) \;.
\end{eqnarray*}
We also set $Z(l,n)=  Z_{l2^{-n} } (2^{-n})$. Then,  we get the following.
\begin{eqnarray*}
x_n & \leq & y_n  +
 \sum_{1\leq l<R_02^n}  \baN \Big( \sum_{^{1\leq i \leq Z(l,n)  } } \!\!\!\! \un_{\{\bm( T_i) \geq \alpha f (2^{-n})  \} }  \Big) \\
& \leq & y_n   +
 \sum_{1\leq l<R_02^n } v(l2^{-n}) \baN_{l2^{-n}}  \Big( \sum_{^{1\leq i \leq Z(l,n)  } } \!\!\!\!
 \un_{\{ \bm(T_i) \geq \alpha f (2^{-n})  \} }  \Big) \\
 & \leq & y_n  +
  \sum_{1\leq  l<R_02^n} v(l2^{-n}) \,  \baN_{l2^{-n}} \big( Z(l, n) \big)  \,
\baN_{2^{-n}} \Big( \bm (B(\rho , 3. 2^{-n} )) \geq  \alpha f(2^{-n})  \Big) \\
& \leq &y_n  +  \sum_{1\leq  l<R_02^n} v(l2^{-n}) \baN_{l2^{-n}} \big( \langle \ell^{l2^{-n}} \rangle  \big)  \, y_n \; .
\end{eqnarray*}
Here, we used (\ref{branchspecify}) in the third line. Recall from (\ref{Nabranchmean}) that $ v(l2^{-n}) \baN_{l2^{-n}} \big( \langle \ell^{l2^{-n}} \rangle  \big) = 1$. Thus, $ x_n \leq  R_0 2^{n} y_n  $. 
We next get an equivalent of $y_n$: 
by (\ref{heightform}) with $\gamma= 2$, we have $v(2^{-n})= 2^n $;  the scaling property (\ref{trianscale}) and Lemma \ref{inverBrown} with $c=2$ imply
\begin{eqnarray*}
y_n & = &  2^n \baN_{2^{-n}} \big( \bm (B(\rho , 3. 2^{-n} )) \geq \alpha f (2^{-n})  \big) \\
&=& 2^n \baN_{1} \big( \bm (B(\rho , 3 )) \geq  \alpha \log 2^{n} \big) \\
& \underset{^{n \rightarrow \infty}}{\sim}  & \frac{_2}{^3} \exp \big( - \big(\frac{_{\pi^2 \log 2}}{^{36}} \alpha -\log 2  \big) \,  n \, \big) \; .
\end{eqnarray*}
Thus,
$$ x_n \leq  R_0 2^{n} y_n  \underset{^{n \rightarrow \infty}}{\sim}  \frac{_2}{^3} R_0
\exp \big( - \big(\frac{_{\pi^2 \log 2}}{^{36}} \alpha -2\log 2 \big) \,  n \, \big) $$
which implies (\ref{WsumBC}) if $\alpha > 72/ \pi^2$.
We argue as in the proof of Theorem \ref{thin} to prove that (\ref{WsumBC}) implies 
\begin{equation}
\label{lowdyadBrown}
\textrm{$\baN$-a.e.}  \qquad \limsup_{n \rightarrow \infty} \; \frac{1}{f (2^{-n})}  \;
\sup_{\sigma \in \cT} \bm \big( B(\sigma , 2^{-n}) \big)  \leq \frac{73}{\pi^2} \; .
\end{equation}

\noi
{\bf Lower Bound.}  Let $R_0$ be a positive integer and let $h_0 \in (0, \infty)$. We also fix $\beta \in (0, \infty)$, that is specified further. We introduce the following event
$$ C_n = \big\{ \Gamma (\cT) > h_0  \big\} \cap \big\{  \sup_{^{\sigma \in B(\rho , R_0)}}
\bm (B (\rho , 2^{-n+3}) )   < \beta  f (2^{-n} ) \,  \big\} $$
We assume that $n \geq 4$, and that $h_0 \geq 2^{-n+1}$. Let $k \geq 1$ be such that
$(k+1)2^{-n} \leq h_0 < (k+2)2^{-n}$. We argue on the event $C_n$: let $\sigma \in \cT$ be such that $d(\rho, \sigma )= h_0$; we apply Lemma \ref{trianinball} with $r= 2^{-n+3}$ (and thus $n_r= n-2$) to prove that there exists
$T \in \cD_{k2^{-n} , 2^{-n} , 2^{-n}} $ such that $T \subset B(\sigma, 2^{-n+3}) $. Thus,
\begin{equation}
\label{AsubABrown}
C_n \subset D_n := \big\{ \Gamma (\cT) > k2^{-n}  \big\} \cap \big\{  \forall T_i \in  \cD_{k2^{-n} , 2^{-n} , 2^{-n}} : \bm (T_i) <  \beta f  (2^{-n} )  \big\} \; .
\end{equation}
To simplify notation we set $Z(k,n)= Z_{k2^{-n}} (2^{-n} ) $. We then use (\ref{branchspecify}) to get
\begin{eqnarray*}
\baN (D_n) & = & \baN (\Gamma (\cT) >k2^{-n} ) \, \baN_{k2^{-n}} \Big( \prod_{^{1\leq i\leq Z(k,n)}}
\un_{\{ \bm (T_i) < \beta f (2^{-n} ) \}}  \Big) \\
& =& v(k2^{-n}) \,  \baN_{k2^{-n}} \left( \,   \baN_{2^{-n}} \! \big( \,  \bm (B(\rho, 2^{-n} ) ) \! < \!  \beta
f (2^{-n}) \,  \big)^{Z(k,n)} \right) \; .
\end{eqnarray*}
Recall that under $\baN_{k2^{-n}}$, conditionally given $\cG_{k2^{-n}}$, the random variable 
$Z(k,n)$ has a Poisson distribution with mean $v(2^{-n}) \langle \ell^{k2^{-n}} \rangle $. We then get
$$\baN (D_n)= v(k2^{-n}) \,  \baN_{k2^{-n}} \big( \exp \big(-z_n \langle \ell^{k2^{-n}} \rangle  \big) \big) \; ,$$
where we have set
$$z_n=
v(2^{-n}) \;  \baN_{2^{-n}} \! \big( \,  \bm (B(\rho, 2^{-n} ) ) \geq   \beta
f (2^{-n}) \,  \big)\; .$$
We next apply  (\ref{heightform}) and (\ref{Nabranch}) with $\gamma= 2$ to get $v(a)= 1/a$ and
$ \baN_a (\exp (-\mu \langle \ell^a \rangle))= (1+a\mu)^{-1}$. 
Since $h_0 \geq 2^{-n+1}$, we have $k2^{-n} \geq h_0/3$ and we get 
\begin{equation}
\label{EnnBestepBrown}
\baN (D_n)= \frac{1}{ k2^{-n} (1+ k2^{-n} z_n)} \leq \frac{1}{\frac{1}{3} h_0 (1+ \frac{1}{3} h_0 z_n)}\; ,
\end{equation}
We next apply
(\ref{kappascale}) and Lemma \ref{inverBrown} with $c=0$, to obtain 
$$ z_n = 2^n \baN_1 \big( \bm (B(\rho, 1)) \geq \beta \log 2^n \big) 
 \underset{^{n \rightarrow \infty}}{\sim} 2 \exp \big( \frac{_{\log 2}}{^4} (4-\pi^2 \beta ) n \big)  \, .$$
Then, for any $\beta < 4/\pi^2$, (\ref{AsubABrown}) and (\ref{EnnBestepBrown}) entail that $  \sum_{n \geq 4} \baN \big( C_n \big) < \infty$.
Thus, for any $h_0, R_0 \in (0, \infty)$, $\baN$-a.e.$\;$for any sufficiently large $n$, $\un_{C_n} = 0$. Since $\baN$-a.e.$ \, 0<\Gamma (\cT) <\infty$, one easily gets
\begin{equation}
\label{uppdyadBrown}
\textrm{$\baN$-a.e.}  \qquad \liminf_{n \rightarrow \infty} \;  \frac{1}{f  (2^{-n})}  \;
\sup_{\sigma \in \cT} \bm \big( B(\sigma , 2^{-n+3}) \big)  \geq \frac{3}{\pi^2}  \; .
\end{equation}
Since $f$ is regularly varying at $0$, (\ref{lowdyadBrown}) and (\ref{uppdyadBrown}) entail Theorem \ref{ThickB}. \cqfd

\subsection{Proof of Proposition \ref{Constant}.}
\label{Constantsec}
Let us fix $\gamma \in (1, 2]$. Recall that
$$ g_\gamma (r)= \frac{r^{\frac{\gamma}{\gamma-1}} }{(\log \log 1/r)^{\frac{1}{\gamma -1}}} \, , \quad
r \in (0, e^{-1}) \; .$$
Recall from (\ref{Mstardef}) the definition of $M^*_r (a)$.
We only need to prove the following: for any $a\in (0, \infty)$,
\begin{equation}
\label{Mstarlower}
\textrm{$\bP$-a.s.} \quad \liminf_{r\rightarrow 0} \frac{M^*_r (a)}{g_\gamma (r)} = \gamma-1 \; .
\end{equation}
Indeed, by (\ref{keymass}), (\ref{Mstarlower}) we get
$$ \baN \left( \int_{\cT} \un_{\{  \liminf_{r\rightarrow 0} \bm (B(\sigma ,r)) /g_\gamma (r) \; \neq \;  \gamma -1
\}}  \bm (d\sigma ) \right) = 0 \; , $$
that immediately entails Proposition \ref{Constant}.

\bigskip

\noi
{\bf Lower bound in (\ref{Mstarlower}).} We fix $a\in (0, \infty)$.
Recall from (\ref{Mstarun}) and (\ref{Mstarscale}) that for any $r\in (0, a]$,
\begin{equation}
\label{scalemmass}
 r^{-\frac{\gamma}{\gamma-1}} M^*_r (a)= M^*_1 (1):= M_* \;
 \end{equation}
and recall Lemma \ref{masszero} that gives the tails of $M_*$ at $0+$. We set
$$ Q = e^{C_\gamma } \sqrt{ \frac{\gamma (\gamma-1)}{2 \pi }  } \; ,  $$
that is the right limit in Lemma \ref{masszero}. Fix $\varrho \in (0,1)$ and $c\in (0, \infty)$. Then, (\ref{scalemmass}) and Lemma \ref{masszero} imply that
\begin{eqnarray*}
\bP \left( \;  M^*_{\varrho^n} (a)  \leq (\gamma-1) c g_\gamma ( \varrho^n) \; \right)  \!\!\!\!\!  & = & \!\!\!\!\!
 \bP \left(    M_* \leq  (\gamma-1) c  (\log \log (\varrho^{-n}))^{-\frac{1}{\gamma -1}} \, \right)  \\
  \!\!\!\!\!  & \underset{^{n \rightarrow \infty}}{\sim} & \!\!\!\!\!
 Q   c^{ \frac{\gamma-1}{2}} (\log 1/\varrho )^{-c^{-(\gamma -1)}}
(\log n)^{-1/2}
\, n^{-c^{-(\gamma -1)}}  .
\end{eqnarray*}
By  Borel-Cantelli, for any $c < 1$, $\bP$-a.s.$\, \liminf_{n \rightarrow \infty}
M^*_{\varrho^n} (a) /g_\gamma (\varrho^n) \geq (\gamma -1) c $. An easy argument, entails that $\bP$-a.s.
$$ \liminf_{n \rightarrow \infty} \frac{M^*_{\varrho^n} (a) }{g_\gamma (\varrho^n)} \geq \gamma -1 \; .$$
For any $r \in (0, 1/\varrho)$, let $n(r) \in \bN$ be such that $\varrho^{n(r)} < r \leq \varrho^{n(r)-1}$. Thus,
$$ \frac{M^*_r(a)  }{g_\gamma (r)}  \geq
\varrho^{\frac{\gamma}{\gamma-1}}
\left(\frac{\log (\log( \varrho^{1-n(r)} ) )}{ \log (\log(
\varrho^{-n(r)} ) )}\right)^{\frac{1}{\gamma -1}}
\frac{ M^*_{ \varrho^{n(r)}} (a)   }{g_\gamma (\varrho^{n(r)})} \; .$$
Thus, for any $\varrho \in (0, 1)$, $\bP$-a.s.$\,  \liminf_{r \rightarrow 0} M^*_r(a) /g_\gamma (r) \geq \varrho^{\frac{\gamma}{\gamma-1}}  (\gamma -1) $, and
by letting $\varrho$ go to $1$, we get
\begin{equation}
\label{lowermass}
\textrm{$\bP$-a.s.} \quad \quad \liminf_{r \rightarrow 0} \; \; \frac{M_r^* (a) }{g_\gamma (r)}
  \geq \gamma -1 \; .
\end{equation}

\noi
{\bf Upper bound in (\ref{Mstarlower}).}
For any $n \geq 2$, we set $r_n= (\log n )^{-n}$ and
$$\varepsilon_n= \un_{\{ M_{ r_n}^* (a)    \leq (\gamma-1)  g_\gamma (r_n)   \}}\; , $$
and $S_n= \varepsilon_2 + \ldots + \varepsilon_n$. Then, (\ref{scalemmass}) and Lemma \ref{masszero} imply that
\begin{equation}
\label{equivepsi}
\bE \left[ \varepsilon_n \right]  \;\underset{^{n \rightarrow \infty}}{\sim} \;
\; Q\,   (\log \log n  )^{-1}
(\log n)^{-1/2} n^{-1} \; .
\end{equation}
Therefore, $\lim_{n \rightarrow \infty} \bE [ S_n] = \infty$. Next observe that
\begin{equation}
\label{carre}
 \bE [ S_n^2] = \bE [ S_n] + 2 \sum_{2 \leq k < l \leq n} \bE [  \varepsilon_k \varepsilon_l] \; .
 \end{equation}
We then use the following lemma.
\begin{lemma}
\label{keymasspac}
There exists a constant $q  \in (0, \infty)$ that only depend
on $\gamma $ such that for any $2 \leq k < l $, $\bE [  \varepsilon_k \varepsilon_l] \leq q   \bE [  \varepsilon_k ] \bE [ \varepsilon_l]$.
\end{lemma}
\noi
{\bf Proof:} first recall that $(U_t, t \geq 0)$ is a subordinator defined on $(\Omega, \cF, \bP)$ with Laplace exponent $\lambda \mapsto \gamma \lambda^{\gamma -1}$. Then, recall
that $\cN^* = \sum_{^{j \in \cI^*}} \delta_{ (  r^*_j , \, H^{*j} )  } $
is a random point measure on $[0, \infty) \times C([0, \infty), \bR)$ defined on  $(\Omega, \cF, \bP)$ such that conditionally given $U$, $\cN^*$ is distributed as a Poisson point measure with intensity $dU_r \otimes \baN (dH)$. Next recall for any $0\leq r^\prime \leq r \leq a$, the notation
$$ M^*_{r^\prime , r} (a) = \sum_{j \in \cI^*} \un_{(r^\prime ,  r]} (r^*_j)  \int_{0}^{\zeta^*_j} \un_{\{ H^{*j}_s \leq r-r^*_j  \}} ds \;  , $$
where $\zeta^*_j$ stands for the lifetime of $H^{*j}$, for any $j \in \cI^*$. Recall from Lemma \ref{basicMprop} that
$M^*_{r^\prime , r} (a) \leq M^*_{ r} (a)$ and that $M^*_{r^\prime , r} (a) $ is independent from
$M^*_{r^\prime } (a) $. Thus, observe that for any $2 \leq k < l $,
\begin{eqnarray*}
 \left\{ \;   M_{r_k}^* (a)    \leq  (\gamma-1) g_\gamma  (r_k)    \; \right\} \cap
 \left\{ \;  M_{r_l}^* (a)   \leq  (\gamma -1) g_\gamma (r_l)    \; \right\}   \hspace{30mm} \\
\hspace{20mm} \subset
\left\{ \;   M^*_{r_l , r_k } (a)   \leq  (\gamma-1) g_\gamma  (r_k)    \; \right\} \cap  \left\{ M_{ r_l}^* (a)   \leq  (\gamma -1) g_\gamma (r_l)    \; \right\} \; .
\end{eqnarray*}
Thus,
$$ \bE [  \varepsilon_k \varepsilon_l] \leq \bP \left( \;    M^*_{r_l , r_k } (a)  \leq
(\gamma-1) g_\gamma  (r_k) \;   \right)   \bE [ \varepsilon_l] \; .$$
Next recall from (\ref{Mstarscale}) that $(r-r^\prime)^{-\frac{\gamma}{\gamma-1}}   M^*_{r^\prime , r} (a) $ has the same law as $M_*$. Consequently,
\begin{eqnarray*}
\bE [  \varepsilon_k \varepsilon_l]  & \leq &  \bP \left(  (1-(r_l/r_k))^{\frac{\gamma}{\gamma -1}} M_*
   \;  \leq \; (\gamma-1) ( \log \log 1/r_k )^{-\frac{1}{\gamma -1}}   \right)   \bE [ \varepsilon_l] \\
& \leq & \bP \left(  (1-(r_{k+1}/r_k))^{\frac{\gamma}{\gamma -1}} M_*    \;  \leq \;
(\gamma-1) ( \log \log 1/r_k )^{-\frac{1}{\gamma -1}}   \right)   \bE [ \varepsilon_l] \; .
\end{eqnarray*}
Now observe that
$$r_{k+1}/r_k=  (\log k)^{-1}  + O ( (\log k)^{-2} ) \quad {\rm and} \quad \log \log1/r_k = \log k + \log \log \log k \; .$$
By Lemma \ref{masszero}, we get
\begin{eqnarray*}  \bP \left( \;   M_*  \leq (\gamma-1) (1-(r_{k+1}/r_k))^{-\frac{\gamma}{\gamma -1}}
 (\log \log 1/r_k )^{-\frac{1}{\gamma -1}}  \; \right)  \; \underset{^{k\rightarrow \infty}}{\sim}
 \hspace{30mm} \\
Q
\big( 1-\frac{_{r_{k+1}}}{^{r_k}} \big)^{-\frac{\gamma}{2}} \left( \log \log 1/r_k  \right)^{-1/2} \exp
\left( - \big( 1-\frac{_{r_{k+1}}}{^{r_k}} \big)^\gamma
 \log\log 1/r_k    \right)  \; ,
\end{eqnarray*}
and since
$$\big( 1-\frac{_{r_{k+1}}}{^{r_k}} \big)^{\gamma}  \log \log 1/r_k  = \log k + \log \log \log k -\gamma +
O \left( \frac{\log \log  \log k  }{\log k}\right) , $$
we get,
\begin{eqnarray*}  \bP \left( \;   M_*  \leq (\gamma-1) (1-(r_{k+1}/r_k))^{-\frac{\gamma}{\gamma -1}}
 (\log \log 1/r_k )^{-\frac{1}{\gamma -1}}  \; \right)  \; \underset{^{k\rightarrow \infty}}{\sim}
 \hspace{30mm} \\
 e^\gamma Q (\log \log k  )^{-1}
(\log k)^{-1/2} k^{-1} \; ,
\end{eqnarray*}
which easily completes the proof of the lemma by (\ref{equivepsi}). \cqfd

\bigskip

The previous lemma and (\ref{carre}) imply that $ \limsup_{n \rightarrow \infty} \bE [ S_n^2 ] /
(  \bE [ S_n] )^2 \leq q  $.
By the Kochen-Stone Lemma, we get $ \bP( \sum_{n \geq 2} \varepsilon_n = \infty ) \geq 1/q  >0$. Namely, with the lower bound (\ref{lowermass}), it entails that
\begin{equation}
\label{quasiMstar}
 \bP \left( \liminf_{r\rightarrow 0} M^*_r (a) / g_\gamma (r) = \gamma -1 \right) \geq 1/ q  >0 \; .
\end{equation}
Standard arguments on Poisson point processes imply that $\liminf_{r\rightarrow 0}M^*_r(a)/g_\gamma (r)$ is measurable with respect to the tail sigma-field of $U$ at $0+$. By Blumenthal zero-one law,
(\ref{quasiMstar}) entails (\ref{Mstarlower}), which completes the proof of Proposition \ref{Constant}. \cqfd

\subsection{Proof of Proposition \ref{constant}.}
\label{constantsec}
Let us fix $\gamma=2$ and let us recall that $g (r)= r^{2} \log \log 1/r $, $ (0, e^{-1}) $. 
Recall from (\ref{Mstardef}) the definition of $M^*_r (a)$. We only need to prove the following: for any $a\in (0, \infty)$,
\begin{equation}
\label{Mstarupper}
\textrm{$\bP$-a.s.} \quad \limsup_{r\rightarrow 0} \frac{M^*_r (a)}{g (r)} = \frac{4}{\pi^2} \; .
\end{equation}
Indeed by the formula (\ref{keymass}), we get
$$ \baN \left( \int_{\cT} \un_{\{  \limsup_{r\rightarrow 0} \bm (B(\sigma ,r)) / g(r) \; \neq \;  4/\pi^2
\}}  \bm (d\sigma ) \right) = 0 \; , $$
that immediately entails Proposition \ref{constant}.

\bigskip

\noi
{\bf Upper bound in (\ref{Mstarupper}).} Recall from (\ref{Mstarscale}) that $r^{-2}M^*_r(a)$ has the law as $M_*$. We fix  $\varrho \in (0, 1)$ and $c \in (0, \infty)$. By Lemma \ref{inverBrownMstar}, we get
$$ \bP ( M_{ \varrho^n}^* (a)  \geq c g (\varrho^n) ) = \bP (M_* \geq  c \log \log \varrho^{-n} ) 
 \underset{^{n \rightarrow \infty}}{\sim} 
4c \, (\log 1/\varrho )^{-\frac{\pi^2}{4}c}   \,  n^{-\frac{\pi^2}{4}c} \log n \; . $$
Borel-Cantelli and an easy argument imply that $\bP$-a.s.
$$ \limsup_{n \rightarrow \infty} \frac{ M^*_{\varrho^n}(a)}{g(\varrho^n)} \leq \frac{4}{\pi^2}  \; .$$
Let $r \in (0, 1)$. There exists $n(r) \in \bN$ such that $\varrho^{n(r)+1} < r \leq \varrho^{n(r)}$. Thus,
$$ \frac{M^*_{ r}(a)}{g(r)} \leq  \frac{1}{\varrho^2} \frac{M_{ \varrho^{n(r)} }^* (a) }{g(\varrho^{n(r)})} \; .$$
Consequently, for any $\varrho \in (0, 1)$, we have $\bP$-a.s.$\, \limsup_{r \rightarrow 0} M^*_r(a)/g(r) \leq 4(\pi \varrho)^{-2}$.
By letting $\varrho $ go to $1$, we get  for any $a\in (0, \infty)$,
\begin{equation}
\label{upperquasi}
\textrm{$\bP$-a.s.} \quad \limsup_{r \rightarrow 0} \frac{M^*_r (a)}{g(r)} \leq \frac{4}{\pi^2} \; .
\end{equation}

\bigskip

\noi
{\bf Lower bound in (\ref{Mstarupper}).} For any $0\leq r^\prime \leq r \leq a$, recall from (\ref{massshelldef})
the definition of $M^*_{r^\prime, r} (a)$ and recall from Lemma \ref{basicMprop}
that for any sequence $(r_n )_{ n \geq 0}$
such that $0\leq r_{n+1} \leq r_n \leq a $, and $\lim_{n \rightarrow \infty} r_n = 0$, the random variables
$M_{r_{n+1} , r_n}^* (a)$ , $n \geq 0$, are independent. Here, we take $r_n = \varrho^n$, with $\varrho \in (0, 1)$. We fix $c \in (0, \infty)$ and for any $n \geq 1$, we set
$$ \varepsilon_n = \un_{\big\{  M^*_{\varrho^{n+1}, \varrho^n} (a) \;  \geq \;  cg (\varrho^n)   \big\}} \; .$$
The scaling property (\ref{Mstarscale}) and Lemma \ref{inverBrownMstar} entail
$$ \bE [ \varepsilon_n]  =  \bP \big( \varrho^{2n} (1-\varrho)^2M_*  \geq c g (\varrho^n)  \, \big) 
 \underset{^{n \rightarrow \infty}}{\sim}  \frac{4c}{(1-\varrho)^2}(\log1/\varrho )^{-
\frac{\pi^2}{4(1-\varrho)^2}c}   \, n^{-\frac{\pi^2}{4(1-\varrho)^2}c}  \log n. $$
If $ c \leq 4(1-\varrho)^2 \pi^{-2}$,  then $ \sum_{n \geq 1} \bE [\varepsilon_n ] = \infty $.
Since the $\varepsilon_n$'s are independent, the usual converse of Borel-Cantelli entails that $\bP$-a.s.$ \,  \sum_{n \geq 1} \varepsilon_n = \infty $.
Now observe that for any $ n\geq 1$,
$$ \varepsilon_n \leq \un_{\big\{  M_{\varrho^n}^*(a ) \geq  c g (\varrho^n) \big\} } \; .$$
Therefore, for any $\varrho \in (0, 1) $ and any $c \in (0, \infty)$ such that $c \leq 4(1-\varrho)^2 \pi^{-2}$
we $\bP$-a.s. have
$$ \limsup_{r\rightarrow \infty} \frac{M^*_r (a) }{g(r)} \geq \limsup_{n \rightarrow \infty} \frac{M_{\varrho^n}^*(a)  }{g(\varrho^n) } \; \geq \; c  \; .  $$
It easily entails the desired lower bound. It completes the proof of (\ref{Mstarupper})
and that of Proposition \ref{constant}. \cqfd

%%%%%
%%%%%
%%%%%
%%%%%
%%%%% BIBLIOGRAPHY.
%%%%%
%%%%%
%%%%%
%%%%%
%%%%%

\end{document}